\documentclass{amsart}
\usepackage{amssymb}
\usepackage{latexsym}
\usepackage{amsthm}
\usepackage{amscd}
\theoremstyle{plain}

\newtheorem{thm}{Theorem}[section]

\newtheorem{lem}[thm]{Lemma}
\newtheorem{prop}[thm]{Proposition}

\newtheorem{defn}{Definition}[section]

\newtheorem{assumpt}[thm]{Assumption}

\numberwithin{equation}{section}

\newcommand{\md}{\mathrm{mod}}
\newcommand{\sgn}{\mathrm{sgn}}

\newcommand{\Col}{\mathrm{Col}}
\newcommand{\Hom}{\mathrm{Hom}}
\newcommand{\tn}{{\tilde{n}}}

\newcommand{\lfl}{\left\lfloor}
\newcommand{\rfl}{\right\rfloor}

\newcommand{\zr}{\mathbf{0}}

\newcommand{\bd}{\mathbf{d}}
\newcommand{\bk}{\mathbf{k}}
\newcommand{\bl}{\mathbf{l}}

\newcommand{\bm}{\mathbf{m}}
\newcommand{\bh}{\mathbf{h}}
\newcommand{\bt}{\mathbf{t}}
\newcommand{\bv}{\mathbf{v}}

\newcommand{\BA}{\mathbf{A}}
\newcommand{\BB}{\mathbf{B}}
\newcommand{\BC}{\mathbf{C}}
\newcommand{\BD}{\mathbf{D}}
\newcommand{\BE}{\mathbf{E}}
\newcommand{\BH}{\mathbf{H}}
\newcommand{\BK}{\mathbf{K}}
\newcommand{\BL}{\mathbf{L}}
\newcommand{\BM}{\mathbf{M}}
\newcommand{\BP}{\mathbf{P}}
\newcommand{\BQ}{\mathbf{Q}}
\newcommand{\BR}{\mathbf{R}}
\newcommand{\BS}{\mathbf{S}}
\newcommand{\BT}{\mathbf{T}}
\newcommand{\BU}{\mathbf{U}}
\newcommand{\BX}{\mathbf{X}}
\newcommand{\BY}{\mathbf{Y}}

\newcommand{\ZZ}{\mathbb Z}
\newcommand{\CC}{\mathbb C}
\newcommand{\QQ}{\mathbb Q}
\newcommand{\RR}{\mathbb R}

\newcommand{\csp}{\null\hskip 20pt}
\newcommand{\ccsp}{\null\hskip 40pt}
\newcommand{\cccsp}{\null\hskip 80pt}
\newcommand{\ccccsp}{\null\hskip 160pt}

\makeatletter
  \def\varddots{\mathinner{\mkern1mu
      \raise\p@\hbox{.}\mkern2mu\raise4\p@\hbox{.}\mkern2mu
      \raise7\p@\vbox{\kern7\p@\hbox{.}}\mkern1mu}}
\makeatother

\newfont{\bg}{cmr10 scaled\magstep2}

\newcommand{\bigzerou}{smash{\lower1.7ex\hbox{\bg 0}}}

\newfont{\bbg}{cmr10 scaled\magstep4}

\newcommand{\bbigzerou}{smash{\lower1.7ex\hbox{\bbg 0}}}

\allowdisplaybreaks[1]

\begin{document}

\title[EXPLICIT FORMULAS FOR HECKE OPERATORS]
{Explicit formulas for Hecke operators on cusp forms,
Dedekind symbols and period polynomials}
\author{Shinji Fukuhara}
\address{Department of Mathematics, Tsuda College, Tsuda-machi 2-1-1,
  Kodaira-shi, Tokyo 187-8577, Japan}
\email{fukuhara@tsuda.ac.jp}
%\date{\today}
\thanks{The author wishes to thank Professor N.~Yui and the referee
for their helpful advice.}

\begin{abstract}
Let $S_{w+2}$ be the vector space of cusp forms
of weight $w+2$
on the full modular group,
and let $S_{w+2}^*$ denote its dual space.
Periods of cusp forms can be regarded as elements of
$S_{w+2}^*$.
The Eichler-Shimura isomorphism theorem asserts that
odd (or even) periods span $S_{w+2}^*$.
However, periods are not linearly independent;
in fact, they satisfy the Eichler-Shimura relations.
This leads to a natural question: which periods would form
a basis of $S_{w+2}^*$.

First we give an answer to this question.
Passing to the dual space $S_{w+2}$,
we will determine a new basis for $S_{w+2}$.
The even period polynomials of this basis elements are expressed
explicitly by means of Bernoulli polynomials.

Next we consider three spaces---$S_{w+2}$,
the space of even Dedekind symbols of weight $w$
with polynomial reciprocity laws,
and the space of even period polynomials of degree $w$.
There are natural correspondences
among these three spaces.
All these spaces are equipped with compatible action of Hecke operators.
We will find explicit form of period polynomials and
the actions of Hecke operators on the period polynomials.

Finally we will obtain
explicit formulas for Hecke operators on $S_{w+2}$
in terms of Bernoulli numbers $B_k$
and divisor functions $\sigma_k(n)$, which are quite
different from the Eichler-Selberg trace formula.
\end{abstract}

\maketitle

\section{Introduction}
\label{sect0}

Throughout the paper,
we assume that $w$ is an {\em even} positive integer,
and we use the following notation:
\begin{align*}
  \Gamma&:=SL_2(\ZZ) \text{\ \ (the full modular group)}, \\
  S_{w+2}&:=
    \text{the vector space of cusp forms on $\Gamma$ with weight $w+2$}, \\
  S_{w+2}^*&:=
    \Hom_\CC( S_{w+2},\CC) \text{\ \ (the dual space of $S_{w+2}$)}, \\
  d_w&:=\begin{cases}
       \lfl\frac{w+2}{12}\rfl-1
         & \mathrm{if\ \ \ } w\equiv 0 \pmod {12} \\
       \lfl\frac{w+2}{12}\rfl
         & \mathrm{if\ \ \ } w\not\equiv 0 \pmod {12}
       \end{cases}
  \end{align*}
where $\lfloor x\rfloor$ denotes the greatest integer not exceeding $x\in\RR$.

One of the main reasons to study cusp forms stems from the fact
that Fourier coefficients of the forms are interesting
from arithmetic view point. Furthermore, periods of
cusp forms (or values of L-series
associated with the forms) play important roles
in number theory.
A large number of papers have discussed periods
and period polynomials of cusp forms
(\cite{E1,F1,F2,H1,KZ1,KZ2,M1,M2,SH1,Z1,Z2}).
One of the basic results is, perhaps, the result due to
Eichler and Shimura (\cite{E1,KZ1,M2,SH1}).

Periods of cusp forms in $S_{w+2}$ can be regarded as elements of
$S_{w+2}^*$.
The Eichler-Shimura isomorphism theorem asserts that
odd (or even) periods span $S_{w+2}^*$ (\cite{E1,KZ1,M2,SH1}).
However, periods are not linearly independent.
In fact, they satisfy the so-called Eichler-Shimura
relations (\cite{LA1,M2}).
This leads us to a natural question: which periods would form
a basis of $S_{w+2}^*$.
The first goal of this paper is to determine odd periods
which constitute a basis of $S_{w+2}^*$ (Theorem \ref{thm1.2}).
Passing to the dual space $S_{w+2}$,
this gives rise to a new basis for $S_{w+2}$.

Next we consider three spaces---$S_{w+2}$,
the space of even Dedekind symbols of weight $w$ with polynomial
reciprocity laws,
and the space of even period polynomials of degree $w$.
It is known that these three spaces are isomorphic
modulo trivial elements (\cite{F2,F4}).
Through these isomorphisms, we can construct bases
for these spaces explicitly
(Theorems \ref{thm1.3}, \ref{thm1.6} and \ref{thm1.7}),
starting with the basis for $S_{w+2}^*$ determined
in Theorem \ref{thm1.2}.
Furthermore, it is also known that these three spaces
are equipped with compatible actions of Hecke operators
(\cite{F4,M2,Z1}).
We will subsequently find explicit forms for the actions
of Hecke operators on the elements of these three spaces
(Theorem \ref{thm1.8}).

As the final goal, we obtain
explicit formulas for Hecke operators on $S_{w+2}$ (Theorem \ref{thm1.9}),
which seem quite different from the celebrated
Eichler-Selberg trace formula.
We will do this by obtaining matrices which represent
the Hecke operators
$T_m\ (m=1,2,\ldots)$ on $S_{w+2}$
as well as their characteristic polynomials
by means of Bernoulli numbers $B_k$
and divisor functions $\sigma_k(n)$.

It might be interesting to compare our result with the
Eichler-Selberg trace formula
(\cite{E1,SE1},\cite[p.\ 48]{LA1}).
Their formula gives traces of Hecke operators in terms of
class numbers of imaginary quadratic fields.
Their method is based on integrating a kernel function for
Hecke operator.
Our approach is different from their method. In fact,
our formulas give
matrices representing
the Hecke operators in terms of Bernoulli numbers and
divisor functions,
as well as their characteristic polynomials.
Our method is based on representing the Hecke operator on $S_{w+2}$
as Hecke operator on the space of Dedekind symbols,
and then as Hecke operator on the space
of period polynomials.
It should be noted that
our argument depends on the fact that $\dim S_{w+2}=d_w$;
while the Eichler-Selberg trace formula gives this fact
as a consequence.

This paper is organized as follows. In Section \ref{sect1},
we give precise statement of our results.
The Sections from \ref{sect2} to \ref{sect5} are devoted to
the preparation of establishing Theorem \ref{thm1.2}.
In Section \ref{sect6}, we will give proofs of Theorems \ref{thm1.2}
and \ref{thm1.3}.
In Section \ref{sect7}, we will present proofs of Theorems \ref{thm1.6}
and \ref{thm1.7}.
The Sections from \ref{sect8} to \ref{sect12} are devoted to
studying the Hecke operators on the three spaces,
which include a proof of Theorem \ref{thm1.8}.
Finally in Section \ref{sect13}, we will prove Theorem \ref{thm1.9}
and append a computer program for obtaining
matrices which represent the Hecke operators.

\section{Statement of results}
\label{sect1}
In this section, we will state our results in more precise form.

Let $f$ be an element of $S_{w+2}$.
We write $f$ as a Fourier series
\begin{equation*}
  f(z)=\sum_{l=1}^{\infty}a_le^{2\pi ilz}.
\end{equation*}
Let $L(f,s)$ be the L-series of $f$.
Namely,
\begin{equation*}
  L(f,s):=\sum_{l=1}^{\infty}\frac{a_l}{l^s}
  \ \ (\Re(s)\gg 0).
\end{equation*}
Then $n$th period of $f$, \ $r_{w,n}(f)$,\  is defined by
\begin{equation*}
  r_{w,n}(f):=\int_{0}^{i\infty}f(z)z^{n}dz
%             =\frac{n!\,i^{n+1}}{(2\pi)^{n+1}}L(f,n+1)\ \ (n=0,1,\ldots,w).
             =\frac{n!}{(-2\pi i)^{n+1}}L(f,n+1)\ \ (n=0,1,\ldots,w).
\end{equation*}
Also the period polynomial of $f$ in the variable $X$ is defined by
\begin{equation*}
  r(f)(X):=\int_{0}^{i\infty}f(z)(X-z)^{w}dz.
\end{equation*}
It is clear that $r(f)(X)$ has the following expression:
\begin{equation*}
  r(f)(X)=\sum_{n=0}^{w}(-1)^n\binom{w}{n}r_{w,w-n}(f)X^{n}.
\end{equation*}
Here and hereafter $\binom{w}{n}$ denotes a binomial coefficient.

Each period $r_{w,n}$ can be regarded as a linear map
from $S_{w+2}$ to $\CC$, that is,
\begin{equation*}
  r_{w,n}\in S_{w+2}^*.
\end{equation*}

Here we recall the result of Eichler and Shimura
(\cite{E1,KZ1,M2,SH1}):
\begin{thm}[Eichler-Shimura]\label{thm1.1}
The maps
\begin{align*}
  r_{w}^+:S_{w+2}&\ \ \to\ \ \CC^{(w+2)/2} \\
     \ f\ \ \ &\ \ \mapsto\ \
     (r_{w,0}(f),\ r_{w,2}(f),\ \ldots,\ r_{w,w}(f))
\end{align*}
and
\begin{align*}
  r_{w}^-:S_{w+2}&\ \ \to\ \ \CC^{w/2} \\
     \ f\ \ \ &\ \ \mapsto\ \
     (r_{w,1}(f),\ r_{w,3}(f),\ \ldots,\ r_{w,w-1}(f))
\end{align*}
are both injective.

In other words,
\begin{enumerate}
\item
the even periods
\begin{equation*}
  r_{w,0},\ r_{w,2},\ \ldots,\ r_{w,w}
\end{equation*}
span the vector space $S_{w+2}^*$;
\item
the odd periods
\begin{equation*}
  r_{w,1},\ r_{w,3},\ \ldots,\ r_{w,w-1}
\end{equation*}
also span $S_{w+2}^*$.
\end{enumerate}
\end{thm}

However, these periods are not linearly independent.
In fact, they satisfy the so-called Eichler-Shimura
relations (\cite{LA1,M2}):
For $n=0,1,\ldots,w$, it holds that
\begin{equation*}\tag{ES1}
  r_{w,n}(f)+(-1)^nr_{w,w-n}(f)
  =0,
\end{equation*}
\begin{equation*}\tag{ES2}
  (-1)^nr_{w,n}(f)
  +\sum_{\substack{0\leq m\leq n \\ m\equiv 0\!\!\!\!\pmod 2}}
    \binom{n}{m}r_{w,w-n+m}(f)
  +\sum_{\substack{0\leq m\leq w-n\\ m\equiv n\!\!\!\!\pmod 2}}
    \binom{w-n}{m}r_{w,m}(f)
  =0,
\end{equation*}
\begin{equation*}\tag{ES3}
  \sum_{\substack{1\leq m\leq n \\ m\equiv 1\!\!\!\!\pmod 2}}
    \binom{n}{m}r_{w,w-n+m}(f)
  +\sum_{\substack{0\leq m\leq w-n\\ m\not\equiv n\!\!\!\!\pmod 2}}
    \binom{w-n}{m}r_{w,m}(f)
  =0.
\end{equation*}

This leads us to a natural question.
Which periods are linearly independent?
Or more strictly, which periods form a basis for $S_{w+2}^*$?
Our first result provides an answer to this question in the case
of odd periods.

Throughout this paper we adopt the following notation and convention:
\begin{defn}\label{defn1.1}
\begin{enumerate}
\item
For an integer $i$ such that $1\leq i\leq d_w$, let $4i\pm 1$
stand for $4i+1$ or $4i-1$
according as $w\equiv 0\pmod 4$\  or \ $w\equiv 2\pmod 4$.
Namely
\begin{equation*}
  4i\pm 1:=\begin{cases}
       4i+1 & \mathrm{if\ \ \ } w\equiv 0\pmod 4 \\
       4i-1 & \mathrm{if\ \ \ } w\equiv 2\pmod 4.
       \end{cases}
\end{equation*}
\item
For an integer $n$ with $0\leq n\leq w$,
let $\tn$ stand for $w-n$:
\begin{equation*}
  \tn:=w-n.
\end{equation*}
\item
For an integer $k$,
a divisor function $\sigma_k$ is defined by
\begin{equation*}
  \sigma_k(n):=\sum_{\substack{ad=n \\ a>0}}a^k,
  \ \ \ (n\in\ZZ^+).
\end{equation*}
\end{enumerate}
\end{defn}
We recall the well-known fact (see e.g. \cite[p.\ 133]{A2}) that
\begin{equation*}
  \dim S_{w+2}=d_w.
\end{equation*}

Now we can state our first result:
\begin{thm}\label{thm1.2}
\begin{equation*}
  \{r_{w,{4i\pm1}}\ |\ i=1,2,\ldots,d_w\}
\end{equation*}
form a basis for $S_{w+2}^*$.
\end{thm}
In other words,
$\{r_{w,{4i\pm1}}\ |\ i=1,2,\ldots,d_w\}$
are linearly independent over $\CC$, and thus other periods are linear
combinations of
$\{r_{w,{4i\pm1}}\ |\ i=1,2,\ldots,d_w\}$.
Furthermore, other odd periods are linear combinations
not only over $\CC$, but over $\QQ$
(confer the proof of Theorem \ref{thm1.2}).

Next we will display a basis for $S_{w+2}$.
For $f,\ g\in S_{w+2}$,
let $(f,g)$ denote
the Petersson inner product.
Then there is a cusp form $R_{w,n}$, which is characterized
by the formula:
\begin{equation*}
  r_{w,n}(f)=(R_{w,n},f)\ \mbox{\ \ for any \ \ } f\in S_{w+2}.
\end{equation*}
Explicit form of $R_{w,n}$,
as a Poincar\'{e} series,
was given (\cite{C1,KZ1}):
\begin{equation*}
  R_{w,n}(z):=
      c_{w,n}^{-1}
      \sum_{\substack{\left[\begin{smallmatrix}a&b\\c&d\end{smallmatrix}\right]
                      \in \Gamma}}
           \frac{1}{(az+b)^{n+1}(cz+d)^{\tn+1}},
      \ \ \ \ c_{w,n}=(-1)^{n+1}2\pi i\binom{w}{n}.
\end{equation*}
In general, the Poincar\'{e} series $R_{w,n}$ is expected to have
transcendental Fourier coefficients.

Passing to the dual space, we obtain a basis of $S_{w+2}$.
\begin{thm}\label{thm1.3}
\begin{equation*}
  \{R_{w,{4i\pm1}}\ |\ i=1,2,\ldots,d_w\}
\end{equation*}
form a basis for $S_{w+2}$.
\end{thm}

Several bases are known for $S_{w+2}$ (\cite{A2,SE2}).
We believe that the above basis is
the first one whose even periods can be described explicitly
(\cite[Theorem $1'$]{KZ1}).
From this fact, we can obtain bases for
the spaces of even period polynomials, as well as
bases for the spaces of even Dedekind symbols
with polynomial reciprocity laws.
We now recall the relationship between
cusp forms, Dedekind symbols and period polynomials.

A complex-valued function $E$ on
$\,\ZZ^+\times\ZZ$
is called a weighted Dedekind symbol of weight $w$
if it satisfies the following two conditions (confer to \cite{F4}):
\begin{equation*}
  E(h,k)=E(h,k+h),\ \ \ E(ch,ck)=c^wE(h,k)
\end{equation*}
for any $(h,k)\in\ZZ^+\times\ZZ$ and $c\in\ZZ^+$.

Moreover, a weighted Dedekind symbol $E$
is said to be {\em even} (resp. {\em odd}) if $E$ satisfies
\begin{equation*}
  E(h,-k)=E(h,k) \text{\ \ \ \ $($resp.\ $E(h,-k)=-E(h,k))$}.
\end{equation*}

There are two rather trivial Dedekind symbols $G_w$ and $F_w$ which are
defined by
\begin{equation*}
  G_w(h,k):=\left\{\gcd(h,k)\right\}^w
  \text{\ \ and\ \ \ }
  F_w(h,k):=h^w
\end{equation*}
for any $(h,k)\in\ZZ^+\times\ZZ$.

A symbol $E$ is determined by its {\em reciprocity law}
\begin{equation*}
  E(h,k)-E(k,-h)=S(h,k)
\end{equation*}
up to addition of scalar multiples of $G_w$.
Here $S$ is a complex-valued function
defined on $\ZZ^+\times\ZZ^+$.

Next we would like to demonstrate the relationship between cusp forms,
weighted Dedekind symbols, and period polynomials.
We need the following notation:
{\allowdisplaybreaks
  \begin{align*}
  \mathcal{W}_w&:=\{
    E\ |\
    \text{$E$ is a Dedekind symbol of weight $w$}
    \}, \\
  \mathcal{W}_w^-&:=\left\{
    E\in\mathcal{W}_w|\ E \text{\ \ is odd\ }
    \right\}, \\
  \mathcal{W}_w^+&:=\left\{
    E\in\mathcal{W}_w|\ E \text{\ \ is even\ }
    \right\}, \\
  \mathcal{E}_w&:=\{
    E\in\mathcal{W}_w|\
    \text{$E(h,k)-E(k,-h)$ is a homogeneous polynomial in $h$ and $k$} \\
    &\ccccsp\cccsp\ccsp\text{of degree $w$}\}, \\
  \mathcal{E}_w^-&:=\left\{
    E\in\mathcal{E}_w|\ E \text{\ \ is odd\ }
    \right\}, \\
  \mathcal{E}_w^+&:=\left\{
    E\in\mathcal{E}_w|\ E \text{\ \ is even\ }
    \right\}, \\
  \mathcal{U}_w&:=\{
    g\ |\
    \text{$g$ is a homogeneous polynomial
          in $h$ and $k$ of degree $w$} \\
    &\ccsp\csp\ \ \text{satisfying $g(h+k,k)+g(h,h+k)=g(h,k)$ and $g(1,1)=0$}
    \} \\
    &\ \ \text{(an element of $\mathcal{U}_w$ is essentially
    a period polynomial modulo $h^w-k^w$ \cite{F2,KZ1}),} \\
  \mathcal{U}_w^-&:=\left\{
    g\in \mathcal{U}_w|\
    \text{\ $g$ is an odd polynomial, i.e., $g(h,-k)=-g(h,k)$}
    \right\}, \\
  \mathcal{U}_w^+&:=\left\{
    g\in \mathcal{U}_w|\
    \text{\ $g$ is an even polynomial, i.e., $g(h,-k)=g(h,k)$}
    \right\}. \\
  \end{align*}
}

For a cusp form $f\in S_{w+2}$ and $(h,k)\in\ZZ^+\times\ZZ$,
we define $E_f$ and $E_f^\pm$ by
\begin{equation}\label{eqn1.1}
  E_f(h,k):=\int_{k/h}^{i\infty}f(z)(hz-k)^{w}dz, \ \ \
  E_f^\pm(h,k):=\frac{1}{2}\{E_f(h,k)\pm E_f(h,-k)\}.
\end{equation}
Then we can show $E_f$ is a Dedekind symbols of weight $w$,
and we can define maps
\begin{equation*}
  \alpha_{w+2}:S_{w+2}\to\mathcal{W}_w,\ \
  \alpha_{w+2}^\pm:S_{w+2}\to\mathcal{W}_w^\pm
\end{equation*}
by
\begin{equation*}
  \alpha_{w+2}(f)=E_f,\ \
  \alpha_{w+2}^\pm(f)=E_f^\pm.
\end{equation*}
Furthermore, we know that $E_f$ and $E_f^{\pm}$ have polynomial
reciprocity laws,
that is, $E_f\in\mathcal{E}_w$ and
$E_f^{\pm}\in\mathcal{E}_w^\pm$.
Hence we have the restricted maps
\begin{equation*}
  \alpha_{w+2}^\pm:S_{w+2}\to\mathcal{E}_w^\pm
\end{equation*}
(to ease the notation, we use the same notation $\alpha_{w+2}^\pm$
for the restricted maps).
Then we have the following:

\begin{thm}[{\cite[Theorem 1.1]{F4}}]\label{thm1.4}
The map
\begin{equation*}
  \alpha_{w+2}^-:S_{w+2}\to\mathcal{E}_w^-
\end{equation*}
is an isomorphism $($between vector spaces$)$, and the map
\begin{equation*}
  \alpha_{w+2}^+:S_{w+2}\to\mathcal{E}_w^+
\end{equation*}
is a monomorphism such that
its image $\alpha_{w+2}^+(S_{w+2})$
is a subspace of $\mathcal{E}_w^+$ of codimension two,
and that $\alpha_{w+2}^+(S_{w+2})$,
$F_w$ and $G_w$ span $\mathcal{E}_w^+$.
\end{thm}

Next we will see how weighted Dedekind symbols are linked
to period polynomials.
For a weighted Dedekind symbol $E$
and $(h,k)\in \ZZ^+\times\ZZ^+$,
let $\beta_{w}(E)$ be defined by
  \begin{equation*}
    \beta_{w}(E)(h,k)=E(h,k)-E(k,-h).
  \end{equation*}
In the case of Dedekind symbol $E_f$ associated with
a cusp form $f$, $\beta_w(E_f)$ has the following expression:
\begin{equation}\label{eqn1.2}
  \beta_w(E_f)(h,k)=\int_{0}^{i\infty}f(z)(hz-k)^{w}dz.
\end{equation}
Note that the right hand side of \eqref{eqn1.2} is
nothing but a homogeneous period polynomial of $f$.

For $E\in\mathcal{E}_w$, we know that $\beta_w(E)\in\mathcal{U}_w$.
Thus, we have a homomorphism
  \begin{equation*}
    \beta_{w}:\mathcal{E}_w\to \mathcal{U}_w.
  \end{equation*}
Then we see that $\beta_{w}$ is almost isomorphism
in the following sense:

\begin{thm}[{\cite[Theorem 1.2]{F4}}]\label{thm1.5}
The homomorphism $\beta_{w}:\mathcal{E}_w\to \mathcal{U}_w$
is an epimorphism with
$\beta_{w}(\mathcal{E}_w^\pm)=\mathcal{U}_w^\pm$,
and $\ker\beta_{w}$ is one dimensional subspace of $\mathcal{E}_w$
spanned by $G_w$.

In particular, the restricted map
  \begin{equation*}
    \beta_{w}^-:\mathcal{E}_w^-\to \mathcal{U}_w^-
  \end{equation*}
is an isomorphism, and
  \begin{equation*}
    \beta_{w}^+:\mathcal{E}_w^+\to \mathcal{U}_w^+
  \end{equation*}
is an epimorphism where
$\ker\beta_{w}^+$ is one dimensional subspace of $\mathcal{E}_w^+$
spanned by $G_w$.
\end{thm}

Here we examine the composed maps
\begin{equation*}
  \beta_{w}^\pm\alpha_{w+2}^\pm:S_{w+2}\to
  \mathcal{E}_w^\pm\to\mathcal{U}_w^\pm.
\end{equation*}
Since $\beta_w^\pm(E_f^\pm)(h,k)=\beta_w^\pm\alpha_{w+2}^\pm(f)(h,k)$
is the homogeneous period polynomial of $f$,
the composed maps
\begin{equation*}
  \beta_{w}^-\alpha_{w+2}^-:S_{w+2}\to
  \mathcal{E}_w^-\to\mathcal{U}_w^-
\end{equation*}
and
\begin{equation*}
  \beta_{w}^+\alpha_{w+2}^+:S_{w+2}\to
  \mathcal{E}_w^+\to\mathcal{U}_w^+
\end{equation*}
can be identified with the Eichler-Shimura
isomorphisms (refer to \cite[p.\ 200]{KZ1}, \cite[Theorem 7.3]{F2}).
In fact,
$\beta_{w}^-\alpha_{w+2}^-$ is an isomorphism, and
$\beta_{w}^+\alpha_{w+2}^+$
is an monomorphism such that the image
$\beta_{w}^+\alpha_{w+2}^+(S_{w+2})$
and $h^w-k^w$ span $\mathcal{U}_w^+$.

These facts may be summarized in the following commutative diagram:

\setlength{\unitlength}{1mm}
\begin{picture}(125,50)(4,-1)
  \put(44,35){\framebox(39,10)
    {\shortstack{the space of cusp forms \\
      of weight $w+2$}}
  }
  \put(0,5){\framebox(65,15)
    {\shortstack{the space of odd (resp. even) Dedekind \\
                 symbols of weight $w$ with polynomial \\
                 reciprocity laws ($\md$ $F_w$ and $G_w$ if even)}}
  }
  \put(77,5){\framebox(50,15)
    {\shortstack{the space of odd (resp. even) \\
                 period polynomials of degree $w$ \\
                 ($\md$ $h^w-k^w$ if even).}}
  }
  \put(44,33){\vector(-1,-1){10}}
  \put(28,27){\makebox(10,5){$\alpha_{w+2}^\pm$}}
  \put(38,25){\makebox(10,5){$\cong$}}
  \put(84,33){\vector(1,-1){10}}
  \put(97,22){\makebox(10,15)
    {\shortstack{the Eichler-Shimura \\
                 isomorphism}}
  }
  \put(80,25){\makebox(10,5){$\cong$}}
  \put(68,12){\vector(1,0){7}}
  \put(66,14){\makebox(10,5){$\beta_{w}^\pm$}}
  \put(66,6){\makebox(10,5){$\cong$}}

  \put(56,-3){\makebox(15,5){Diagram ES}}
\end{picture}
\vspace*{1mm}

Using these correspondences,
we can obtain bases for the spaces
of even period polynomials,
and bases for the spaces of even Dedekind symbols with polynomial
reciprocity laws.
We need the following notation (refer to \cite{F3}).

For an integer $n$ such that $0<n<w$, a polynomial $S_{w,n}$,
in $h$ and $k$, is defined by
{\allowdisplaybreaks
\begin{align*}
    S_{w,n}(h,k):=&
       (-1)^n\frac{B_{n+1}(\frac{k}{h})h^w}{n+1}
        +\frac{B_{n+1}(\frac{h}{k})k^w}{n+1}
        -\frac{B_{\tn+1}(\frac{k}{h})h^w}{\tn+1}
        -(-1)^n\frac{B_{\tn+1}(\frac{h}{k})k^w}{\tn+1} \\
      &\ccsp\ \ \ \ \ \ +
      \begin{cases}
      \frac{w+2}{B_{w+2}}
      \frac{B_{n+1}}{n+1}\frac{B_{\tn+1}}{\tn+1}(h^w-k^w)
         & \mathrm{if\ \ \ } n\equiv 1 \pmod {2} \\
      0  & \mathrm{if\ \ \ } n\equiv 0 \pmod {2}.
      \end{cases}
\end{align*}
}

For any $(h,k)\in \ZZ^+\times\ZZ$,
and for an integer $n$ such that $0<n<w$,
a Dedekind symbol $E_{w,n}:\ZZ^+\times\ZZ\to \CC$ is defined by
{\allowdisplaybreaks
\begin{align*}
  E_{w,n}(h,k)
    :=&\frac{1}{2}
      \sum_{\substack{\left[\begin{smallmatrix}a&b\\c&d\end{smallmatrix}\right]
                      \in \Gamma \\
                      ac\ne 0 \\
                      (k/h+b/a)(k/h+d/c)<0}}
      \sgn\left(\frac{k}{h}+\frac{b}{a}\right)
        (ak+bh)^{\tn}(ck+dh)^{n} \\
      &\cccsp+
      \left\{
        (-1)^n\frac{\bar{B}_{n+1}(\frac{k}{h})h^w}{n+1}
        -\frac{\bar{B}_{\tn+1}(\frac{k}{h})h^w}{\tn+1}
        \right\} \\
      &\cccsp+
      \begin{cases}
      \frac{w+2}{B_{w+2}}
      \frac{B_{n+1}}{n+1}\frac{B_{\tn+1}}{\tn+1}h^w
         & \mathrm{if\ \ \ } n\equiv 1 \pmod {2} \\
      0  & \mathrm{if\ \ \ } n\equiv 0 \pmod {2}.
      \end{cases}
\end{align*}
}
Here and hereafter, $B_m(x)$ (resp. $B_m$) denotes
the $m$th Bernoulli polynomial (resp. number), and
$\bar{B}_m(x)$ denotes the $m$th Bernoulli function.
That is, $B_m(x)$ is defined by
\begin{equation*}
  \frac{te^{xt}}{e^t-1}=\sum_{m=0}^{\infty}B_m(x)\frac{t^m}{m!},
\end{equation*}
and $\bar{B}_m(x)$ is defined as the periodic function which
coincides with $B_m(x)$ on $[0,1)$.
Moreover, $\sgn(x)$ denotes the sign of $x\in\RR$.

It is shown (\cite[Theorem $1'$]{KZ1}, \cite{F2}) that,
for $n$ odd$, \ \beta_{w}^+\alpha_{w+2}^+$
maps $c_{w,n}R_{w,n}$ to $S_{w,n}$:
\begin{equation*}
  \beta_{w}^+\alpha_{w+2}^+ : c_{w,n}R_{w,n}\mapsto S_{w,n}
    \ \ (n\mbox{\ odd}).
\end{equation*}
Furthermore, we will show in Lemma \ref{lem7.1} below that,
for $n$ odd$, \ \alpha_{w+2}^+$ maps $c_{w,n}R_{w,n}$ to $E_{w,n}$:
\begin{equation*}
  \alpha_{w+2}^+ : c_{w,n}R_{w,n}\mapsto E_{w,n} \ \ (n\mbox{\ odd}).
\end{equation*}

From these facts and Theorems \ref{thm1.3}, \ref{thm1.4}
and \ref{thm1.5}, we obtain the following
two theorems:
\begin{thm}\label{thm1.6}
\begin{equation*}
  \{\ S_{w,{4i\pm1}}(h,k)\ |\ i=1,2,\ldots,d_w\ \}\cup\{\ h^w-k^w\ \}
\end{equation*}
form a basis for $\mathcal{U}_w^+$.
\end{thm}

\begin{thm}\label{thm1.7}
\begin{equation*}
  \{\ E_{w,{4i\pm1}}\ |\ i=1,2,\ldots,d_w\ \}\cup\{\ F_w\ \}\cup\{\ G_w\ \}
\end{equation*}
form a basis for $\mathcal{E}_w^+$.
\end{thm}

The latter half of this paper is devoted to obtaining
matrices which represent the Hecke operators $T_m\ (m=1,2,\ldots)$
on $S_{w+2}$ as well as their characteristic polynomials.
For this purpose, first we discuss Hecke operators
on the three spaces in Diagram ES.
Manin \cite{M2} (see also \cite[p.\ 202]{KZ1}) and Zagier \cite{Z1} proved
that there are well-defined Hecke operators (also denoted by $T_m$)
on the spaces of period polynomials which are compatible with
the Eichler-Shimura isomorphism:

\begin{equation}\label{eqn1.3}
  \begin{CD}
  S_{w+2} @>{\beta_{w}^\pm\alpha_{w+2}^\pm}>> \mathcal{U}_w^\pm \\
  @VV{T_m}V @VV{T_m}V \\
  S_{w+2} @>{\beta_{w}^\pm\alpha_{w+2}^\pm}>> \mathcal{U}_w^\pm.
  \end{CD}
\end{equation}

Furthermore, in \cite{F4}, we introduced
Hecke operators on the space of
Dedekind symbols by the following formula:
\begin{equation*}
  (T_{m}E)(h,k)
    :=\sum_{\substack{ad=m \\ d>0}}
      \sum_{b(\md d)}E(dh,ak+bh).
\end{equation*}

It was proved that
Hecke operators on the spaces of Dedekind symbols are compatible
with Hecke operators on the spaces of cusp forms (\cite{F4}).
As a consequence,
we have the following commutative diagram:

\begin{equation}\label{eqn1.4}
  \begin{CD}
  S_{w+2} @>{\alpha_{w+2}^\pm}>> \mathcal{E}_w^\pm \\
  @VV{T_m}V @VV{T_m}V \\
  S_{w+2} @>{\alpha_{w+2}^\pm}>> \mathcal{E}_w^\pm.
  \end{CD}
\end{equation}

(To ease the notation, we will use the same notation $T_m$
for the Hecke operators on $S_{w+2}$, $\mathcal{E}_w^\pm$
and $\mathcal{U}_w^\pm$.)

Now we need the following definitions
to describe the actions of Hecke operators on
$R_{w,n}$, $E_{w,n}$ and $S_{w,n}$:
\begin{defn}\label{defn1.2}
\begin{enumerate}
\item
For a positive integer $m$,
\begin{equation*}
  H_m:=\left\{
       \begin{bmatrix}a&b\\c&d\end{bmatrix}
       \ |\ \ ad-bc=m;\ a,b,c,d\in\ZZ
       \right\};
\end{equation*}
\item
For positive integers $m$ and $n$ such that $0<n<w$,
\begin{equation*}
  R_{w,n}^m(z):=
      m^{w+1}c_{w,n}^{-1}
      \sum_{\substack{\left[\begin{smallmatrix}a&b\\c&d\end{smallmatrix}\right]
                      \in H_m}}
           \frac{1}{(az+b)^{n+1}(cz+d)^{\tn+1}};
\end{equation*}
\item
For positive integers $m$ and $n$ such that $0<n<w$,
we define a map $E_{w,n}^m:\ZZ^+\times\ZZ\to \CC$
by
{\allowdisplaybreaks
\begin{align*}
  E_{w,n}^m(h,k)
    :=&\frac{1}{2}
      \sum_{\substack{\left[\begin{smallmatrix}a&b\\c&d\end{smallmatrix}\right]
                      \in H_m \\
                      ac\ne 0 \\
                      (k/h+b/a)(k/h+d/c)<0}}
      \sgn\left(\frac{k}{h}+\frac{b}{a}\right)
        (ak+bh)^{\tn}(ck+dh)^{n} \\
      &\ccsp+
      \sum_{\substack{ad=m \\ a>0}}
      \left\{(-1)^n
        d^{\tn}\frac{\bar{B}_{n+1}(\frac{ak}{h})h^w}{n+1}
        -d^n\frac{\bar{B}_{\tn+1}(\frac{ak}{h})h^w}{\tn+1}
        \right\} \\
      &\ccsp+
      \begin{cases}
      \sigma_{w+1}(m)\frac{w+2}{B_{w+2}}
      \frac{B_{n+1}}{n+1}\frac{B_{\tn+1}}{\tn+1}h^w
         & \mathrm{if\ \ \ } n\equiv 1 \pmod {2} \\
      0  & \mathrm{if\ \ \ } n\equiv 0 \pmod {2};
      \end{cases}
\end{align*}
}
\item
For positive integers $m$ and $n$ such that $0<n<w$,
we define a polynomial $S_{w,n}^m$
in $h$ and $k$ by
{\allowdisplaybreaks
\begin{align*}
    S_{w,n}^m(h,&k)
    :=\frac{1}{2}
      \sum_{\substack{\left[\begin{smallmatrix}a&b\\c&d\end{smallmatrix}\right]
                      \in H_m \\
                      abcd<0}}
      \sgn(ab)(ak+bh)^{\tn}(ck+dh)^{n} \\
      &\ccsp+\sum_{\substack{ad=m \\ a>0}}
        d^\tn\left\{
        (-1)^n\frac{B_{n+1}(\frac{ak}{h})h^w}{n+1}
        +\frac{B_{n+1}(\frac{ah}{k})k^w}{n+1}
        \right\} \\
      &\ccsp-\sum_{\substack{ad=m \\ a>0}}
        d^n\left\{
        \frac{B_{\tn+1}(\frac{ak}{h})h^w}{\tn+1}
        +(-1)^n\frac{B_{\tn+1}(\frac{ah}{k})k^w}{\tn+1}
        \right\} \\
      &\ccsp
      +\begin{cases}
      \sigma_{w+1}(m)\frac{w+2}{B_{w+2}}
      \frac{B_{n+1}}{n+1}\frac{B_{\tn+1}}{\tn+1}(h^w-k^w)
         & \mathrm{if\ \ \ } n\equiv 1 \pmod {2} \\
      0  & \mathrm{if\ \ \ } n\equiv 0 \pmod {2}.
      \end{cases}
\end{align*}
}
\end{enumerate}
\end{defn}

Properties of these functions will be studied in the later sections.
It is plain that
\begin{equation*}
  R_{w,n}^1=R_{w,n},\ \ E_{w,n}^1=E_{w,n}\ \
  \text{and\ \ }S_{w,n}^1=S_{w,n}.
\end{equation*}

Now we can formulate the actions of Hecke operators on
$R_{w,n}$, $E_{w,n}$ and $S_{w,n}$:
\begin{thm}\label{thm1.8}
The actions of the Hecke operators $T_m$ on
$R_{w,n}$, $E_{w,n}$ and $S_{w,n}$
are expressed as follows:
\begin{equation*}
  T_m(R_{w,n})=R_{w,n}^m,\ \
  T_m(E_{w,n})=E_{w,n}^m\
  \mathrm{\ \ and\ \ \ }
  T_m(S_{w,n})=S_{w,n}^m.
\end{equation*}
\end{thm}

Finally, as an application of Theorems \ref{thm1.3} and \ref{thm1.8},
we will give explicit formulas for the Hecke operators
on the spaces of cusp forms.
Let
\begin{equation*}
  f(h,k)=\sum_{\nu=0}^{w}a_\nu h^\nu k^{w-\nu}
  \text{\ \ \ and\ \ \ }
  g(h,k)=\sum_{\nu=0}^{w}b_\nu h^\nu k^{w-\nu}
\end{equation*}
be homogeneous polynomials in $h$ and $k$ with degree $w$.
Then their inner product
$\langle f,g\rangle$ is defined by
\begin{equation*}
  \langle f,g\rangle:=\sum_{\nu=0}^{w}a_\nu \bar{b}_\nu
\end{equation*}
where $\bar{b}_\nu$ denotes the complex conjugate of $b_\nu$.

Under this notation we obtain the following result.
\begin{thm}\label{thm1.9}
\begin{enumerate}
\item
Let $m$ be a positive integer, and let $\BT_{m}$ be
the matrix representing the Hecke operator
\begin{equation*}
  T_m:S_{w+2}\to S_{w+2}
\end{equation*}
with respect to the basis
\begin{equation*}
  c_{w,4i\pm 1}R_{w,{4i\pm 1}}\ \ (i=1,2,\ldots,d_w).
\end{equation*}
Let $\BS_1$ and $\BS_2$ be matrices defined by
\begin{equation*}
  \BS_1:=
   \begin{bmatrix}
   \langle S_{w,{4i\pm 1}},S_{w,{4j\pm 1}}\rangle
   \end{bmatrix}
   \ \ \ (i,j=1,2,\ldots,d_w),
\end{equation*}
\begin{equation*}
  \BS_2:=
   \begin{bmatrix}
   \langle S_{w,{4i\pm 1}},S_{w,{4j\pm 1}}^m\rangle
   \end{bmatrix}
   \ \ \ (i,j=1,2,\ldots,d_w).
\end{equation*}
Then $\BT_{m}$ can be expressed as
\begin{equation*}
  \BT_{m}=\BS_1^{-1}\BS_2.
\end{equation*}
\item
Let $n$ be an odd integer with $0<n<w$. Then $S_{w,{n}}^m$ can be
expressed explicitly as a polynomial in $h$ and $k$ by the following
formula:
{\allowdisplaybreaks
\begin{align*}
  S_{w,{n}}^m(h,k)
    =&
    2\sum_{\substack{\nu=0 \\ \nu\mathrm{\ even} }}^w
    \Biggl\{
    \sum_{\mu=1}^{m-1}
    \sum_{\lambda=\max(0,\nu-\tn)}^{\min(n,\nu)}
    \mu^{\lambda}(\mu-m)^{n-\lambda}
    \binom{\tn}{\nu-\lambda}\binom{n}{\lambda}\times \\
    &\cccsp\csp\ \ \  \sigma_{\tn-\nu}(\mu)\sigma_{\nu-n}(m-\mu)
    \Biggr\}
    h^{\nu}k^{w-\nu} \\
      &+\frac{(-1)^nm^\tn}{n+1}\sum_{\nu=\tn-1}^{w}\binom{n+1}{\nu-\tn+1}
        B_{\nu-\tn+1}\sigma_{n-\nu}(m)h^{\nu}k^{w-\nu} \\
      &-\frac{m^n}{\tn+1}\sum_{\nu=n-1}^{w}\binom{\tn+1}{\nu-n+1}
        B_{\nu-n+1}\sigma_{\tn-\nu}(m)h^{\nu}k^{w-\nu} \\
      &+\frac{m^\tn}{n+1}\sum_{\nu=0}^{n+1}\binom{n+1}{n-\nu+1}
        B_{n-\nu+1}\sigma_{\nu-\tn}(m)h^{\nu}k^{w-\nu} \\
      &-\frac{(-1)^nm^n}{\tn+1}\sum_{\nu=0}^{\tn+1}\binom{\tn+1}{\tn-\nu+1}
        B_{\tn-\nu+1}\sigma_{\nu-n}(m)h^{\nu}k^{w-\nu} \\
      &+\sigma_{w+1}(m)\frac{w+2}{B_{w+2}}
      \frac{B_{n+1}}{n+1}\frac{B_{\tn+1}}{\tn+1}(h^w-k^w).
\end{align*}
}
In particular, setting $m=1$, we have
{\allowdisplaybreaks
\begin{align*}
  S_{w,{n}}(h,k)
      =&\frac{(-1)^n}{n+1}\sum_{\nu=\tn-1}^{w}\binom{n+1}{\nu-\tn+1}
        B_{\nu-\tn+1}h^{\nu}k^{w-\nu} \\
      &-\frac{1}{\tn+1}\sum_{\nu=n-1}^{w}\binom{\tn+1}{\nu-n+1}
        B_{\nu-n+1}h^{\nu}k^{w-\nu} \\
      &+\frac{1}{n+1}\sum_{\nu=0}^{n+1}\binom{n+1}{n-\nu+1}
        B_{n-\nu+1}h^{\nu}k^{w-\nu} \\
      &-\frac{(-1)^n}{\tn+1}\sum_{\nu=0}^{\tn+1}\binom{\tn+1}{\tn-\nu+1}
        B_{\tn-\nu+1}h^{\nu}k^{w-\nu} \\
      &+\frac{w+2}{B_{w+2}}
      \frac{B_{n+1}}{n+1}\frac{B_{\tn+1}}{\tn+1}(h^w-k^w).
\end{align*}
}
\end{enumerate}
\end{thm}

Consequently, this shows that we can express the matrix
$\BT_{m}$ representing the Hecke operator $T_{m}$ explicitly
in terms of Bernoulli numbers $B_k$
and divisor functions $\sigma_k(n)$.

In the last section, we append a computer program
for obtaining matrices which represent
the Hecke operators,
and their characteristic polynomials.
The program is a straightforward implementation of Theorem \ref{thm1.9}.

\section{The Eichler-Shimura relations for modified periods}
\label{sect2}

The Sections from \ref{sect2} to \ref{sect6}
are devoted to the study of odd periods of cusp forms.
We also present proofs of Theorems \ref{thm1.2} and \ref{thm1.3}
in Section \ref{sect6}.

Our strategy for proving Theorems \ref{thm1.2} is as follows.
We obtain an integral matrix which express
the Eichler-Shimura relations for odd periods. Then we take
the reduction modulo 2 of this matrix. The new matrix has a nice
``self-similar'' structure so that we can find linearly independent
column vectors. Then we choose odd periods corresponding to the
column vectors which turn out to form a basis.

For our purpose, it is more convenient
to consider ``modified periods'' instead of periods themselves.
\begin{defn}\label{defn2.1}
For $f\in S_{w+2}$, we define $s_{w,n}(f)$ by
\begin{equation}\label{eqn2.1}
  s_{w,n}(f):=(-1)^n\binom{w}{n}r_{w,w-n}(f),
\end{equation}
and we call $s_{w,n}(f)$ the $n$th modified period of $f$.
\end{defn}
We regard $s_{w,n}$ as an element of $S_{w+2}^*$, that is,
\begin{equation*}
  s_{w,n}\in S_{w+2}^*.
\end{equation*}
Then the period polynomial $r(f)(X)$ has the following expression:
\begin{equation*}
  r(f)(X)=\sum_{n=0}^{w}(-1)^n\binom{w}{n}r_{w,w-n}(f)X^{n}
         =\sum_{n=0}^{w}s_{w,n}(f)X^{n}.
\end{equation*}

Using the modified period $s_{w,n}$,
the Eichler-Shimura relations can be expressed
as follows (Kohnen-Zagier \cite[p.\ 199]{KZ1}):
\begin{equation*}\tag{KZ1}
  s_{w,n}+(-1)^ns_{w,w-n}
  =0\ \ \ \ (0\leq n\leq w),
\end{equation*}
\begin{equation*}\tag{KZ2}
  \sum_{\substack{m=0 \\ m\ \mathrm{even}}}^{n}\binom{w-m}{w-n}s_{w,m}
  +\sum_{\substack{m=n+1 \\ m\ \mathrm{even}}}^{w}\binom{m}{n}s_{w,m}
  =0\ \ \ \ (0\leq n\leq w),
\end{equation*}
\begin{equation*}\tag{KZ3}
  \sum_{\substack{m=0 \\ m\ \mathrm{odd}}}^{n}\binom{w-m}{w-n}s_{w,m}
  +\sum_{\substack{m=n+1 \\ m\ \mathrm{odd}}}^{w}\binom{m}{n}s_{w,m}
  =0\ \ \ \ (0\leq n\leq w).
\end{equation*}
These linear relations are the starting point of our discussions
which eventually lead to Theorem \ref{thm1.2}.

Hereafter we adopt the following notation and convention:
\begin{enumerate}
\item
\begin{equation*}
  w_2:=\frac{w}{2} \mbox{\ \ \ and\ \ \ }
  w_4:=\lfl\frac{w}{4}\rfl;
\end{equation*}
\item
Let $\BX=[x_{ij}]$ be an $m\times n$ matrix.
Then, to the end of Section \ref{sect6},
we employ the convention that
the index $i$ (resp. $j$) runs from $0$ to $m-1$ (resp. from $0$ to $n-1$).
\end{enumerate}
We also need the following notation:
\begin{defn}\label{defn2.2}
We define $t_{w,i}\in S_{w+2}^*\ (i=0,1,\ldots,w_2)$ by
\begin{equation}\label{eqn2.2}
  t_{w,i}:=
    \begin{cases}
      0          & \mathrm{if\ \ \ } i=0 \\
      s_{w,2i-1} & \mathrm{if\ \ \ } 1\leq i\leq w_2.
    \end{cases}
\end{equation}
\end{defn}
So, $t_{w,i}\ (i=1,2,\ldots,w_2)$ are odd modified periods,
while $t_{w,0}$ is a dummy introduced for technical reasons.

Our first task is to express the relations (KZ1) and
(KZ3) for odd modified period in matrix forms.
For this, we first introduce the following matrices.

\begin{defn}\label{defn2.3}
\begin{enumerate}
\item
We define a $(w_2+1)\times 1$ matrix $\bt$ by
\begin{equation*}
  \bt:=\begin{bmatrix}
      t_{w,0} \\ t_{w,1} \\ t_{w,2} \\ \vdots \\ t_{w,w_2}
      \end{bmatrix};
\end{equation*}
\item
We define a matrix $\BA=[a_{ij}]$
$($the matrix whose $i$th row and $j$th column entry is $a_{ij}$;
\ $i=0,1,\ldots,w;\ j=0,1,\ldots,w_2$$)$
by
\begin{equation*}
  a_{ij}:=\begin{cases}
      0 & \mathrm{if\ \ \ } j=0 \\
      \binom{2j-1}{i} & \mathrm{if\ \ \ } j\ne0,\ 2j>i \\
      \binom{w-2j+1}{w-i} & \mathrm{if\ \ \ } j\ne0,\ 2j\leq i;
      \end{cases}
\end{equation*}
\item
We define a matrix $\BB=[b_{ij}]\ (i=0,1,\ldots,w_2;\ j=0,1,\ldots,w_2)$ by
\begin{equation*}
  b_{ij}:=\begin{cases}
      0 & \mathrm{if\ \ \ } j=0 \\
      1
      & \mathrm{if\ \ \ } j\ne0,\ i=0 \\
      \binom{2j-1}{2i}+\binom{2j-1}{2i-1}
      & \mathrm{if\ \ \ } j\ne0,\ i\ne 0,\ j>i \\
      \binom{w-2j+1}{w-2i}+\binom{w-2j+1}{w-2i+1}
      & \mathrm{if\ \ \ } j\ne0,\ i\ne 0,\ j\leq i.
      \end{cases}
\end{equation*}
\end{enumerate}
\end{defn}

Then the Eichler-Shimura relations (KZ3) are expressed as
\begin{equation*}
  \BA\bt=\zr.
\end{equation*}
Here and hereafter the symbol $\zr$ stands for a zero matrix.

Furthermore we have:
\begin{lem}\label{lem2.1}
It holds that
\begin{equation*}
  \BB\bt=\zr.
\end{equation*}
\end{lem}
\begin{proof}
We note that $0$th row of $\BB$ is equal to $0$th row of $\BA$,
and that the $i$th row of $\BB$ is
the sum of the $(2i-1)$th row and $(2i)$th
row of $\BA$ for $1\leq i\leq w_2$.
Thus, from the fact that $\BA\bt=\zr$, we have that $\BB\bt=\zr$.
\end{proof}

Next we need the following definitions:
\begin{defn}\label{defn2.4}
We define a matrix $\BC=[c_{ij}]$
$(i=0,1,\ldots,w_4;\ j=0,1,\ldots,w_2)$
as follows,
which depends on the congruence conditions:
\begin{enumerate}
\item
For an integer $i$ with $0\leq i<w_4$, we define
\begin{equation*}
  c_{ij}:=\begin{cases}
      0 & \mathrm{if\ \ \ } j=0 \\
      1 & \mathrm{if\ \ \ } j=i+1 \\
      -1 & \mathrm{if\ \ \ } j=w_2-i \\
      0 & \mathrm{otherwise};
      \end{cases}
\end{equation*}
\item
For $w\equiv 0\pmod 4$ and $i=w_4$, we define
\begin{equation*}
  c_{ij}:=\begin{cases}
      0 & \mathrm{if\ \ \ } 0\leq j\leq w_4 \\
      1 & \mathrm{otherwise};
      \end{cases}
\end{equation*}
\item
For $w\equiv 2\pmod 4$ and $i=w_4$, we define
\begin{equation*}
  c_{ij}:=\begin{cases}
      0 & \mathrm{if\ \ \ } 0\leq j\leq w_4 \\
      1 & \mathrm{if\ \ \ } j=w_4+1 \\
      2 & \mathrm{otherwise}.
      \end{cases}
\end{equation*}
\end{enumerate}
\end{defn}

Here are two examples of the matrix $\BC$:
\begin{enumerate}
\item
For $w=8$,
\begin{equation*}
  \BC=\begin{bmatrix}
      0 & 1 & 0 & 0 & -1 \\
      0 & 0 & 1 & -1 & 0 \\
      0 & 0 & 0 & 1 & 1
      \end{bmatrix};
\end{equation*}
\item
For $w=10$,
\begin{equation*}
  \BC=\begin{bmatrix}
      0 & 1 & 0 & 0 & 0 & -1 \\
      0 & 0 & 1 & 0 & -1 & 0 \\
      0 & 0 & 0 & 1 & 2  & 2 \\
      \end{bmatrix}.
\end{equation*}
\end{enumerate}

From the Eichler-Shimura relations (KZ1) and (KZ3) for $n=0$,
we obtain:
\begin{lem}\label{lem2.2}
\begin{equation}\label{eqn2.3}
  \BC\bt=\zr.
\end{equation}
\end{lem}

Finally, we introduce the following matrix $\BD$,
which can be obtained from the matrices $\BB$ and $\BC$:
\begin{defn}\label{defn2.5}
\begin{equation*}
  \BD:=\begin{bmatrix}
     \BB \\ \BC
     \end{bmatrix}.
\end{equation*}
More precisely,
for $\BD=[d_{ij}]\ (i=0,1,\ldots,w_2+w_4+1;\ j=0,1,\ldots,w_2)$,
\begin{equation*}
  d_{ij}:=\begin{cases}
          b_{ij} & \mathrm{if\ \ \ } 0\leq i\leq w_2,\ 0\leq j\leq w_2 \\
          c_{(i-w_2-1)j}
            & \mathrm{if\ \ \ } w_2+1\leq i\leq w_2+w_4+1,\ 0\leq j\leq w_2.
     \end{cases}
\end{equation*}
\end{defn}

Then, from lemmas \ref{lem2.1} and \ref{lem2.2}, we have
\begin{lem}\label{lem2.3}
\begin{equation*}
  \BD\bt=\zr.
\end{equation*}
\end{lem}

\section{The Eichler-Shimura relations modulo 2}
\label{sect3}

In this section, we consider the reductions modulo 2 of the
matrices which were introduced in the previous section.
By $\ZZ/2\ZZ$, we denote the set of congruence classes modulo 2.
For an integer $x$,
\begin{equation*}
  \bar{x} \text{\ \ \ \ \ \ or\ \ \ \ \ \ } x\!\!\mod 2
\end{equation*}
denotes the congruence class of $x$ modulo $2$
so that we have
\begin{equation*}
  \ZZ/2\ZZ=\{\bar{0},\bar{1}\}.
\end{equation*}

Let $\BK$, $\BL$ and $\BM$ denote the reductions modulo 2
of $\BB$, $\BC$ and $\BD$, respectively:
\begin{equation*}
  \BK:=[k_{ij}]=[\bar{b}_{ij}],\ \
  \BL:=[\ell_{ij}]=[\bar{c}_{ij}]
  \mathrm{\ \ \ and\ \ \ }
  \BM:=[m_{ij}]=[\bar{d}_{ij}].
\end{equation*}
From Definition \ref{defn2.5}, we know that
\begin{equation*}
  \BM=\begin{bmatrix}
     \BK \\ \BL
     \end{bmatrix}.
\end{equation*}

Here we recall Lucas' congruence theorem on binomial coefficients
(\cite[p.\ 271]{DI1}).
Let $p$ be a prime number, and let $n,k,a,b$ be nonnegative integers
with $0\leq a,b<p$. Then it holds that
\begin{equation}\label{eqn3.1}
  \binom{np+a}{kp+b}\equiv\binom{n}{k}\binom{a}{b}\ \ \pmod p.
\end{equation}
Using this identity we can prove the following lemma:
\begin{lem}\label{lem3.1}
As for the matrix $\BK=[k_{ij}]$,
it holds that
\begin{equation}\label{eqn3.2}
  k_{ij}=\begin{cases}
         \binom{j}{i} \mod 2 & \mathrm{if\ \ \ } j>i \\
         \bar{0} & \mathrm{if\ \ \ } j\leq i.
         \end{cases}
\end{equation}
\end{lem}
\begin{proof}
If $j>i>0$, we have
\begin{align*}
  b_{ij}&=\binom{2j-1}{2i}+\binom{2j-1}{2i-1}
        =\binom{2(j-1)+1}{2i}+\binom{2(j-1)+1}{2(i-1)+1} \\
        &\equiv\binom{j-1}{i}\binom{1}{0}
         +\binom{j-1}{i-1}\binom{1}{1}
         \pmod 2 \text{\ \ \ (by \eqref{eqn3.1})} \\
        &=\binom{j-1}{i}+\binom{j-1}{i-1}
        =\binom{j}{i}
         \text{\ \ \ \ \ (by\  Pascal's\  identity).} \\
\end{align*}
This implies that, if $j>i>0$, then
\begin{equation*}
  k_{ij}=\binom{j}{i} \mod 2.
\end{equation*}
If \ $0<j\leq i$, we have
\begin{align*}
  b_{ij}&=\binom{w-2j+1}{w-2i}+\binom{w-2j+1}{w-2i+1} \\
        &\equiv\binom{w_2-j}{w_2-i}\binom{1}{0}
         +\binom{w_2-j}{w_2-i}\binom{1}{1}
         \pmod 2 \text{\ \ \ (by \eqref{eqn3.1})} \\
        &\equiv 0 \pmod 2.
\end{align*}
This shows that, if $\ 0<j\leq i$, then
\begin{equation*}
  k_{ij}=\bar{0}.
\end{equation*}

In the case that $i=0$ or $j=0$, the identities follow
from the definitions of $b_{ij}$ and $k_{ij}$.
This completes the proof.
\end{proof}

\section{The Pascal-Sierpinski's triangle and the matrix $\BK$}
\label{sect4}

Here we investigate the matrix $\BK$.
Since $\BK=[k_{ij}]$ satisfies
\begin{equation*}
  k_{ij}=\begin{cases}
         \binom{j}{i} \mod 2 & \mathrm{if\ \ \ } j>i \\
         \bar{0} & \mathrm{if\ \ \ } j\leq i,
         \end{cases}
\end{equation*}
$\BK$ is an upper triangular matrix,
and the upper triangular part of the matrix is nothing but
the Pascal-Sierpinski's triangle (see e.g. \cite{B1}).
In particular it is ``self-similar''.
The fact that the Pascal-Sierpinski's triangle has ``self-similarity''
have been well-known (the reader can refer to,
e.g. \cite[pp.\ 44,53]{ST1}
for the congruence properties of binomial coefficients, which gives rise
to the ``self-similarity'').

Considering these, we inductively define a family of square matrices
$\BP_n\ (n=0,1,\ldots)$,
with entries in $\ZZ/2\ZZ$ as follows:
\begin{defn}\label{defn4.1}
\begin{enumerate}
\item
\begin{equation*}
\BP_0:=
\begin{bmatrix}
    \bar{1}
\end{bmatrix};
\end{equation*}
\item
For any positive integer $n$,
\begin{equation*}
\BP_{n}:=
\begin{bmatrix}
    \BP_{n-1}  & \BP_{n-1} \\
    \zr     & \BP_{n-1}
\end{bmatrix}.
\end{equation*}

Note that the size of $\BP_n$ is $2^n\times 2^n$.
\item
For any positive integer $n$,
let $\BE_n$ denote the identity matrix of size $n\times n$.
Namely
\begin{equation*}
  \BE_n:=[\bar{\delta}_{ij}]\ \ (i,j=0,1,\ldots,n-1)
\end{equation*}
where $\delta_{ij}$ is the Kronecker's delta.
\item
For any positive integer $n$,
\begin{equation*}
\BQ_{n}:=\BP_{n}+\BE_{2^n}.
\end{equation*}

Note that the size of $\BQ_n$ is again $2^n\times 2^n$.
\item
Let $k$ and $n$ be positive integers with $k\leq n$, and let
$\BX=[x_{ij}]\ \ (i,j=0,1,\ldots,n-1)$ be $n\times n$-matrix. Then,
by $\BX[[k]]$,
we denote the principal submatrix of $\BX$ with size $k\times k$.
Namely
\begin{equation*}
\BX[[k]]=[x_{ij}]\ \ (i,j=0,1,\ldots,k-1).
\end{equation*}
\end{enumerate}
\end{defn}

Under the notation we can express $\BK$ as follows:
\begin{lem}\label{lem4.1}
Let $n$ be an integer such that $w_2+1\leq 2^n$. Then $\BK$ can be expressed as
\begin{equation*}
  \BK=\BQ_n[[w_2+1]].
\end{equation*}
\end{lem}
\begin{proof}
We note that the size of $\BK$ is $(w_2+1)\times(w_2+1)$
while the size of $\BQ_n$ is $2^n\times 2^n$.
Then the lemma follows from Lemma \ref{lem3.1} which asserts that
\begin{equation*}
  k_{ij}=\begin{cases}
         \binom{j}{i} \mod 2 & \mathrm{if\ \ \ } j>i \\
         \bar{0} & \mathrm{if\ \ \ } j\leq i.
         \end{cases}
\end{equation*}
\end{proof}

In the rest of this section, we prove several lemmas
which will be used in the proof of Theorem \ref{thm1.2}.
The following is obvious:
\begin{lem}\label{lem4.2}
For any positive integer $n$,
\begin{equation*}
\BQ_{n}=
\begin{bmatrix}
    \BQ_{n-1}  & \BQ_{n-1}+\BE_{2^{n-1}} \\
    \zr        & \BQ_{n-1}
\end{bmatrix}.
\end{equation*}
\end{lem}

Here we introduce a special type of operation on matrices.
\begin{defn}\label{defn4.2}
Let $\BX=[x_{ij}], \ \BX'=[x'_{ij}]\ \
(i=0,1,\ldots,n-1; j=0,1,\ldots,k-1)$ be
two matrices of size $n\times k$.
Then $\BX'$ is said to be obtained from $\BX$ by an $R^+$-operation
if there exist integers $i_0$ and $i_1$ with $0\leq i_0<i_1<n$
such that
\begin{equation*}
  x'_{ij}=\begin{cases}
          x_{ij} & \mathrm{if\ \ \ } i\ne i_0 \\
          x_{i_0j}+x_{i_1j} & \mathrm{if\ \ \ } i=i_0.
         \end{cases}
\end{equation*}
In other words, $\BX'$ is obtained from $\BX$ by
adding $i_1$th row to $i_0$th row
with $i_0<i_1$
$($this is an elementary row operation$)$.
\end{defn}
In what follow, the notation
\begin{equation*}
  \BX\Rightarrow \BY
\end{equation*}
means ``$\BY$ is obtained from $\BX$ by a sequence of $R^+$-operations''.
It is clear that $\BX\Rightarrow \BY$ if and only if
$\BY$ can be expressed as
\begin{equation*}
  \BY=
  \begin{bmatrix}
    \bar{1}    & *      & \cdots  & *  \\
    \bar{0}      & \bar{1}    & \ddots  & \vdots       \\
    \vdots & \ddots & \ddots  & *      \\
    \bar{0}      & \cdots & \bar{0}       & \bar{1}
  \end{bmatrix}
  \BX.
\end{equation*}

Then we can formulate the following lemma:
\begin{lem}\label{lem4.3}
For any positive integer $n$, it holds that
\begin{equation*}
  \begin{bmatrix}
    \BE_{2^{n+1}} \\
    \BQ_{n+1}
  \end{bmatrix}
  \Rightarrow
  \begin{bmatrix}
    \BE_{2^n} & \zr       \\
    \BQ_n     & \zr       \\
    \BQ_n     & \BE_{2^n} \\
    \zr       & \BQ_n
  \end{bmatrix}.
\end{equation*}
\end{lem}
\begin{proof}
Using Lemma \ref{lem4.2}, we have
\begin{equation*}
  \begin{bmatrix}
    \BE_{2^n} & \zr       & \zr       & \zr       \\
    \zr       & \BE_{2^n} & \BE_{2^n} & \BE_{2^n} \\
    \zr       & \zr       & \BE_{2^n} & \BE_{2^n} \\
    \zr       & \zr       & \zr       & \BE_{2^n} \\
  \end{bmatrix}
  \begin{bmatrix}
    \BE_{2^{n+1}} \\
    \BQ_{n+1}
  \end{bmatrix}
  =
  \begin{bmatrix}
    \BE_{2^n} & \zr       \\
    \BQ_n     & \zr       \\
    \BQ_n     & \BE_{2^n} \\
    \zr       & \BQ_n
  \end{bmatrix}.
\end{equation*}
\end{proof}

Now we set
\begin{equation*}
  \BU:=
  \begin{bmatrix}
    \BE_2 \\
    \BQ_1
  \end{bmatrix}
     =
  \begin{bmatrix}
    \bar{1} & \bar{0} \\
    \bar{0} & \bar{1} \\
    \bar{0} & \bar{1} \\
    \bar{0} & \bar{0}
  \end{bmatrix}.
\end{equation*}
Then we obtain:
\begin{lem}\label{lem4.4}
For any nonnegative integer $n$, the matrix
\begin{equation*}
  \begin{bmatrix}
    \BE_{2^{n+1}} \\
    \BQ_{n+1}
  \end{bmatrix}
\end{equation*}
can be transformed to a block matrix of the form
\begin{equation*}
  \begin{bmatrix}
    \BU    & \zr      & \cdots  & \zr  \\
    *      & \BU    & \ddots  & \vdots       \\
    \vdots & \ddots & \ddots  & \zr      \\
    *      & \cdots & *       & \BU
  \end{bmatrix}
\end{equation*}
by a sequence of $R^+$-operations.
\end{lem}
\begin{proof}
Repeatedly applying Lemma \ref{lem4.3}, we can transform the matrix
\begin{equation*}
  \begin{bmatrix}
    \BE_{2^{n+1}} \\
    \BQ_{n+1}
  \end{bmatrix}
\end{equation*}
by $R^+$-operations:
\begin{align*}
  \begin{bmatrix}
    \BE_{2^{n+1}} \\
    \BQ_{n+1}
  \end{bmatrix}
  \ \ &\Rightarrow\ \
  \begin{bmatrix}
    \BE_{2^n} & \zr       \\
    \BQ_n     & \zr       \\
    *       & \BE_{2^n} \\
    *       & \BQ_n
  \end{bmatrix}
  \ \ \Rightarrow\ \
  \begin{bmatrix}
    \BE_{2^{n-1}} & \zr           & \zr           & \zr           \\
    \BQ_{n-1}     & \zr           & \zr           & \zr           \\
    *           & \BE_{2^{n-1}} & \zr           & \zr           \\
    *           & \BQ_{n-1}     & \zr           & \zr           \\
    *           & *           & \BE_{2^{n-1}} & \zr           \\
    *           & *           & \BQ_{n-1}     & \zr           \\
    *           & *           & *           & \BE_{2^{n-1}} \\
    *           & *           & *           & \BQ_{n-1}
  \end{bmatrix} \\
  \ \ &\Rightarrow\ \
  \cdots
  \ \ \Rightarrow\ \
  \begin{bmatrix}
    \BE_2    & \zr      & \cdots  & \cdots  & \zr  \\
    \BQ_1    & \zr      & \cdots  & \cdots  & \zr  \\
    *      & \BE_2    & \zr       &         & \vdots \\
    *      & \BQ_1    & \zr       & \ddots  & \vdots \\
    \vdots & \ddots & \ddots  & \ddots  & \zr      \\
    \vdots &        & \ddots  & \ddots  & \zr      \\
    *      & \cdots & *       & *       & \BE_2    \\
    *      & \cdots & *       & *       & \BQ_1
  \end{bmatrix}
  =
  \begin{bmatrix}
    \BU    & \zr      & \cdots  & \zr  \\
    *      & \BU    & \ddots  & \vdots       \\
    \vdots & \ddots & \ddots  & \zr      \\
    *      & \cdots & *       & \BU
  \end{bmatrix}.
\end{align*}
This completes the proof.
\end{proof}

From Lemma \ref{lem4.4}, we can obtain the following:
\begin{lem}\label{lem4.5}
Let $n$ be a positive integer. Then,
by a sequence of $R^+$-operations, $\BQ_{n+2}$ can be transformed
to a block matrix $\check{\BQ}_{n+2}$ of the form
\begin{equation}\label{eqn4.1}
  \check{\BQ}_{n+2}=
  \begin{bmatrix}
  \BX        &
  \begin{bmatrix}
    \BU    & \zr      & \cdots  & \zr  \\
    *      & \BU    & \ddots  & \vdots       \\
    \vdots & \ddots & \ddots  & \zr      \\
    *      & \cdots & *       & \BU
  \end{bmatrix}
          &
  \begin{bmatrix}
  \bar{1}       & *      & \cdots & *      \\
  \bar{0}       & \bar{1}      & \ddots & \vdots \\
  \vdots  & \ddots & \ddots & *      \\
  \bar{0}       & \cdots & \bar{0}      & \bar{1}      \\
  \end{bmatrix} \\
  \zr       &
  \zr       &     \BQ_{n+1}
  \end{bmatrix},
\end{equation}
where $\BX$ is a $2^{n+1}\times 2^n$ matrix.
Furthermore we can assume that $R^+$-operations used in this transformation
are elementary row operations adding $i$th rows with $i<2^{n+1}$.
\end{lem}
\begin{proof}
First we see that $\BQ_{n+2}$ have the form
\begin{align*}
  \BQ_{n+2}&=
  \begin{bmatrix}
    \BQ_{n+1} & \BQ_{n+1}+\BE_{2^{n+1}} \\
    \zr     & \BQ_{n+1}         \\
  \end{bmatrix} \\
  &=
  \begin{bmatrix}
  \begin{bmatrix}
    \BQ_{n} & \BQ_{n}+\BE_{2^n} \\
    \zr     & \BQ_{n}         \\
  \end{bmatrix}
          &
  \begin{bmatrix}
  \bar{1}       & *      & \cdots & *      \\
  \bar{0}       & \bar{1}      & \ddots & \vdots \\
  \vdots  & \ddots & \ddots & *      \\
  \bar{0}       & \cdots & \bar{0}      & \bar{1}      \\
  \end{bmatrix} \\
  \zr       &  \BQ_{n+1}
  \end{bmatrix}.
\end{align*}

Now, applying $R^+$-operations,
this matrix can be transformed to a matrix of the form
\begin{equation*}
  \begin{bmatrix}
  \begin{bmatrix}
    \BQ_{n} & \BE_{2^n}       \\
    \zr     & \BQ_{n}         \\
  \end{bmatrix}
          &
  \begin{bmatrix}
  \bar{1}       & *      & \cdots & *      \\
  \bar{0}       & \bar{1}      & \ddots & \vdots \\
  \vdots  & \ddots & \ddots & *      \\
  \bar{0}       & \cdots & \bar{0}      & \bar{1}      \\
  \end{bmatrix} \\
  \zr       &  \BQ_{n+1}
  \end{bmatrix}
  =
  \begin{bmatrix}
  \begin{bmatrix}
    \BQ_{n}        \\
    \zr            \\
  \end{bmatrix}
          &
  \begin{bmatrix}
    \BE_{2^n}       \\
    \BQ_{n}         \\
  \end{bmatrix}
          &
  \begin{bmatrix}
  \bar{1}       & *      & \cdots & *      \\
  \bar{0}       & \bar{1}      & \ddots & \vdots \\
  \vdots  & \ddots & \ddots & *      \\
  \bar{0}       & \cdots & \bar{0}      & \bar{1}      \\
  \end{bmatrix} \\
  \zr       & \zr      &  \BQ_{n+1}
  \end{bmatrix}.
\end{equation*}

Furthermore, by Lemma \ref{lem4.4}, this matrix can be transformed to
a matrix, say $\BR$, of the form
\begin{equation*}
  \BR=
  \begin{bmatrix}
  \BX        &
  \begin{bmatrix}
    \BU    & \zr      & \cdots  & \zr  \\
    *      & \BU    & \ddots  & \vdots       \\
    \vdots & \ddots & \ddots  & \zr      \\
    *      & \cdots & *       & \BU
  \end{bmatrix}
          &
  \begin{bmatrix}
  \bar{1}       & *      & \cdots & *      \\
  \bar{0}       & \bar{1}      & \ddots & \vdots \\
  \vdots  & \ddots & \ddots & *      \\
  \bar{0}       & \cdots & \bar{0}      & \bar{1}      \\
  \end{bmatrix} \\
  \zr       &
  \zr       &     \BQ_{n+1}
  \end{bmatrix},
\end{equation*}
where $\BX$ is a matrix of size $2^{n+1}\times 2^n$.
We take this matrix $\BR$ as $\check{\BQ}_{n+2}$.
Here we used the fact that, while $R^+$-operations, a matrix of the form
\begin{equation*}
  \begin{bmatrix}
  \bar{1}       & *      & \cdots & *      \\
  \bar{0}       & \bar{1}      & \ddots & \vdots \\
  \vdots  & \ddots & \ddots & *      \\
  \bar{0}       & \cdots & \bar{0}      & \bar{1}      \\
  \end{bmatrix} \\
\end{equation*}
remains a matrix of  this form.
We note that, in this transformation,
we only used elementary row operations adding $i$th rows with $i<2^{n+1}$.

This completes the proof.
\end{proof}

Moreover, from Lemma \ref{lem4.5}, we can obtain the following:
\begin{lem}\label{lem4.6}
Let $n$ be a positive integer. Then,
by a sequence of $R^+$-operations, $\BQ_{n+3}$ can be transformed
to a block matrix $\Check{\Check\BQ}_{n+3}$ of the form
\begin{equation}\label{eqn4.2}
  \Check{\Check\BQ}_{n+3}=
  \begin{bmatrix}
  \begin{bmatrix}
  \BX        &
  \begin{bmatrix}
    \BU    & \zr      & \cdots  & \zr  \\
    *      & \BU    & \ddots  & \vdots       \\
    \vdots & \ddots & \ddots  & \zr      \\
    *      & \cdots & *       & \BU
  \end{bmatrix}
          &
  \begin{bmatrix}
  \bar{1}       & *      & \cdots & *      \\
  \bar{0}       & \bar{1}      & \ddots & \vdots \\
  \vdots  & \ddots & \ddots & *      \\
  \bar{0}       & \cdots & \bar{0}      & \bar{1}      \\
  \end{bmatrix} \\
  \zr       &
  \zr       &     \BQ_{n+1}
  \end{bmatrix}
    & \BY \\
  \zr & \BQ_{n+2}
  \end{bmatrix}
\end{equation}
where $\BX$ and $\BY$ are  $2^{n+1}\times 2^n$ and
$2^{n+2}\times 2^{n+2}$ matrices, respectively.
Furthermore we can assume that $R^+$-operations used in this transformation
are elementary row operations adding $i$th rows with $i<2^{n+1}$.
\end{lem}
\begin{proof}
From Lemmas \ref{lem4.2} and \ref{lem4.5}, we have
\begin{align*}
  \BQ_{n+3}&=
  \begin{bmatrix}
    \BQ_{n+2} & \BQ_{n+2}+\BE_{2^{n+2}} \\
    \zr     & \BQ_{n+2}         \\
  \end{bmatrix}
  \Rightarrow
  \begin{bmatrix}
    \check{\BQ}_{n+2} & \BY \\
    \zr     & \BQ_{n+2}         \\
  \end{bmatrix}. \\
\end{align*}

This completes the proof.
\end{proof}

Hereafter we assume the following:
\begin{assumpt}\label{assumpt4.7}
\begin{equation}\label{eqn4.3}
  1< n
  \mathrm{\ \ \ and\ \ \ }
  2^{n}\leq \frac{w_2}{3} <2^{n+1}.
\end{equation}
\end{assumpt}
Note that, for a given $w_2\geq 12$, there exists uniquely an integer $n$
which satisfies \eqref{eqn4.3}.
Under Assumption \ref{assumpt4.7}, we see obviously that
$w_2+1\leq 2^{n+3}$. Hence we know that $\BK$ can be express
as
\begin{equation*}
  \BK=\BQ_{n+3}[[w_2+1]].
\end{equation*}
Now we set
\begin{equation}\label{4.3}
  \BH:=\Check{\Check\BQ}_{n+3}[[w_2+1]].
\end{equation}

Then, from Lemma \ref{lem4.6}, we obtain
\begin{lem}\label{lem4.8}
By a sequence of $R^+$-operations, the matrix $\BK$ can be
transformed to the matrix $\BH$.
\end{lem}
\begin{proof}
By Lemma \ref{lem4.6}, $\BQ_{n+3}$ can be transformed
to the matrix $\Check{\Check\BQ}_{n+3}$
by a sequence of $R^+$-operations adding $i$th rows with $i<2^{n+1}$.
Furthermore, from Assumption \ref{assumpt4.7}, it follows that
$2^{n+1}\leq w_2+1$. Hence we have proved that
the matrix $\BK$ can be transformed to the matrix $\BH$
by a sequence of $R^+$-operations.
\end{proof}

\section{Properties of the matrix $\BH$}
\label{sect5}

In this section, keeping Assumption \ref{assumpt4.7},
we will study properties of the matrix
\begin{equation*}
  \BH=\Check{\Check\BQ}_{n+3}[[w_2+1]]
\end{equation*}
which is obtained from $\BK$ by a sequence of $R^+$-operations.

We need the following notation:
\begin{defn}\label{defn5.1}
\begin{enumerate}
\item
For an integer $i\geq 0$, we define an integer $\alpha(i)$ by
\begin{equation*}
  \alpha(i):=\lfl\frac{i}{2}\rfl+
  \begin{cases}
      1  & \mathrm{if\ \ \ } i\equiv 1 \pmod 4 \\
      0  & \mathrm{otherwise};
  \end{cases}
\end{equation*}
\item
Let us set
\begin{equation*}
  w=12k+2a\ \ (0<k,\ 0\leq a\leq 5).
\end{equation*}
We define
\begin{equation*}
  \beta_1(w):=-2^{n+1}+
  \begin{cases}
      4k+2  & \mathrm{if\ \ \ } a=0,1 \\
      4k+4  & \mathrm{if\ \ \ } a=2,3,4 \\
      4k+6  & \mathrm{if\ \ \ } a=5;
  \end{cases}
\end{equation*}
\begin{equation*}
  \beta_2(w):=-2^{n+1}+
  \begin{cases}
      4k    & \mathrm{if\ \ \ } a=0,1,2 \\
      4k+2  & \mathrm{if\ \ \ } a=3,4,5.
  \end{cases}
\end{equation*}
\end{enumerate}
\end{defn}

Obviously we have that
\begin{equation*}
  \alpha(i+4)=\alpha(i)+2, \ \ \
  \alpha(i)=\frac{i}{2}
  \mathrm{\ \ if\ \ }
  i\equiv 0\pmod 2; \\
\end{equation*}
\begin{equation}\label{lem5.2-1}
  2^{n}+\alpha(\beta_1(w))=
  \begin{cases}
      2k+1  & \mathrm{if\ \ \ } a=0,1 \\
      2k+2  & \mathrm{if\ \ \ } a=2,3,4 \\
      2k+3  & \mathrm{if\ \ \ } a=5;
  \end{cases} \\
\end{equation}
\begin{equation}\label{lem5.2-2}
  2^{n}+\alpha(\beta_2(w))=
  \begin{cases}
      2k   & \mathrm{if\ \ \ } a=0,1,2 \\
      2k+1 & \mathrm{if\ \ \ } a=3,4,5.
  \end{cases}
\end{equation}

Using the notation $\alpha(i)$, we can formulate
the following lemma:
\begin{lem}\label{lem5.1}
Let us set
\begin{equation*}
\BH=[h_{ij}]\ \ (i,j=0,1,\ldots,w_2),
\end{equation*}
and let $i$ be an integer such that
\begin{equation*}
  0\leq i< 2^{n+1}.
\end{equation*}
Then it holds that
\begin{align}
  h_{ij}&=\bar{0} &\mathrm{for\ \ \ }
         &2^{n}+\alpha(i)<j<\min(\,2^{n+1}+i,w_2+1), \label{eqn5.1}\\
  h_{ij}&=\bar{1} &\mathrm{for\ \ \ }
         &j=2^{n+1}+i,\ \ j< w_2+1, \label{eqn5.2}\\
  h_{ij}&=\bar{1} &\mathrm{for\ \ \ }
         &j=2^n+\alpha(i),\ \ i\equiv 0\pmod 2.\label{eqn5.3}
\end{align}
\end{lem}
\begin{proof}
Let us set
\begin{equation*}
  \Check{\Check\BQ}_{n+3}=[\Check{\Check q}_{ij}]
  \ \ (i,j=0,1,\ldots,2^{n+3}-1).
\end{equation*}
Then, for $i$ such that $0\leq i< 2^{n+1}$,
we can easily read off the following identities
\begin{align}
  \Check{\Check q}_{ij}&=\bar{0} &\mathrm{for\ \ \ }
         &2^{n}+\alpha(i)<j<2^{n+1}+i,\label{eqn5.4}\\
  \Check{\Check q}_{ij}&=\bar{1} &\mathrm{for\ \ \ }
         &j=2^{n+1}+i,\label{eqn5.5}\\
  \Check{\Check q}_{ij}&=\bar{1} &\mathrm{for\ \ \ }
         &j=2^n+\alpha(i),\ \ i\not\equiv 3\pmod 4 \label{eqn5.6}
\end{align}
from the form
\begin{equation*}
  \Check{\Check\BQ}_{n+3}=
  \begin{bmatrix}
  \begin{bmatrix}
      \overbrace{\hspace*{0.5cm}}^{2^n}
            & \overbrace{\hspace{3.0cm}}^{2^n}
            & \overbrace{\hspace{2.5cm}}^{2^{n+1}} \\
  \BX        &
  \begin{bmatrix}
    \BU    & \zr      & \cdots  & \zr  \\
    *      & \BU    & \ddots  & \vdots       \\
    \vdots & \ddots & \ddots  & \zr      \\
    *      & \cdots & *       & \BU
  \end{bmatrix}
          &
  \begin{bmatrix}
  \bar{1}       & *      & \cdots & *      \\
  \bar{0}       & \bar{1}      & \ddots & \vdots \\
  \vdots  & \ddots & \ddots & *      \\
  \bar{0}       & \cdots & \bar{0}      & \bar{1}      \\
  \end{bmatrix} \\
  \zr       &
  \zr       &     \BQ_{n+1}
  \end{bmatrix}
    & \BY \\
  \zr & \BQ_{n+2}
  \end{bmatrix}
\end{equation*}
and
\begin{equation*}
  \BU=
  \begin{bmatrix}
    \bar{1} & \bar{0} \\
    \bar{0} & \bar{1} \\
    \bar{0} & \bar{1} \\
    \bar{0} & \bar{0}
  \end{bmatrix}
\end{equation*}
where $\BX$ (resp. $\BY$) is  a $2^{n+1}\times 2^n$
(resp. $2^{n+2}\times 2^{n+2}$) matrix.

Next, noting that $\BH=\Check{\Check\BQ}_{n+3}[[w_2+1]]$, we have
the identities
\eqref{eqn5.1}, \eqref{eqn5.2} and \eqref{eqn5.3}
from \eqref{eqn5.4}, \eqref{eqn5.5} and \eqref{eqn5.6}
respectively.
This completes the proof.
\end{proof}

By Lemma \ref{lem5.1}, we obtain:

\begin{lem}\label{lem5.3}
Set
$w=12k+2a\ \ (0<k,\ 0\leq a\leq 5)$, and let
\begin{equation*}
  \BH=[h_{ij}]\ \ (i,j=0,1,\ldots,w_2).
\end{equation*}

\begin{enumerate}
\item
Suppose that
\begin{equation*}
  w<2\cdot 2^{n+2}
  \mathrm{\ \ \ and\ \ \ }
  i=\beta_1(w).
\end{equation*}
Then we have
\begin{align}
  h_{ij}&=\bar{0} &\mathrm{if\ \ \ }
  &\begin{cases}
      2k+1<j<4k+2, &\ a=0,1 \\
      2k+2<j<4k+4, &\ a=2,3,4 \\
      2k+3<j<4k+6, &\ a=5,
  \end{cases} \label{eqn5.7}\\
  h_{ij}&=\bar{1} &\mathrm{if\ \ \ }
  &\begin{cases}
      j=2k+1, &\ a=0,1 \\
      j=2k+2, &\ a=2,3,4 \\
      j=2k+3, &\ a=5.
  \end{cases} \label{eqn5.8}
\end{align}
\item
Suppose that
\begin{equation*}
  2\cdot 2^{n+2}\leq w
  \mathrm{\ \ \ and\ \ \ }
  i=\beta_2(w).
\end{equation*}
Then we have
\begin{align}
  h_{ij}&=\bar{0} &\mathrm{if\ \ \ }
  &\begin{cases}
      2k<j<4k, &\ a=0,1,2 \\
      2k+1<j<4k+2, &\ a=3,4,5,
  \end{cases} \label{eqn5.9}\\
  h_{ij}&=\bar{1} &\mathrm{if\ \ \ }
  &\begin{cases}
      j=4k, &\ a=0,1,2 \\
      j=4k+2, &\ a=3,4,5.
  \end{cases} \label{eqn5.10}
\end{align}
\end{enumerate}
\end{lem}
\begin{proof}
Under Assumption \ref{assumpt4.7},
we can easily show the following inequalities
(the detail is left to the reader):
\begin{equation}\label{eqn5.11}
  0\leq\beta_1(w)<2^{n+1}
  \mathrm{\ \ \ if\ \ \ }
  w<2\cdot 2^{n+2},
\end{equation}
\begin{equation}\label{eqn5.12}
  0\leq\beta_2(w)<2^{n+1}
  \mathrm{\ \ \ if\ \ \ }
  2\cdot 2^{n+2}\leq w.
\end{equation}

Furthermore, when $i=\beta_1(w)$ or $i=\beta_2(w)$,
we can easily check that
\begin{equation*}
  2^{n+1}+i< w_2+1.
\end{equation*}

Now we are ready to apply Lemma \ref{lem5.1}
to prove the lemma.
Noting that $\beta_1(w)$ is even,
we can apply \eqref{eqn5.1} and \eqref{eqn5.3} in
Lemma \ref{lem5.1} and
\eqref{lem5.2-1} to obtain \eqref{eqn5.7} and \eqref{eqn5.8}.

We can also apply \eqref{eqn5.1} and \eqref{eqn5.2} in
Lemma \ref{lem5.1} and
\eqref{lem5.2-2} to obtain \eqref{eqn5.9} and \eqref{eqn5.10}.
This completes the proof.
\end{proof}

In the sequel, we use the following notation;
for column vectors $\bv_1,\ldots,\bv_k$, let
\begin{equation*}
  \Col(\bv_1,\ldots,\bv_k):= \text{the column vector space
     spanned by $\bv_1,\ldots,\bv_k$}.
\end{equation*}
We also denote
the $i$th entry of $\bv_j$ by $\bv_j(i)\ \ (i=0,1,\ldots)$.
Let $\bh_j$, $\bk_j$, $\bl_j$ and $\bm_j$ denote column vectors of
$\BH$, $\BK$, $\BL$ and $\BM$, that is
\begin{equation*}
  \BH=[\bh_0\ \ldots\ \bh_{w_2}],\
  \BK=[\bk_0\ \ldots\ \bk_{w_2}],\
  \BL=[\bl_0\ \ldots\ \bl_{w_2}] \text{\ \ and\ \ }
  \BM=[\bm_0\ \ldots\ \bm_{w_2}].
\end{equation*}
Then, for $\bh_j$ and $\bk_j$,
it holds that $\bh_j(i)=h_{ij}$ and $\bk_j(i)=k_{ij}$, respectively.

Under this notation, we have the following lemma:
\begin{lem}\label{lem5.4}
$\bm_1,\bm_2,\ldots,\bm_{w_4}$ are linearly independent.
\end{lem}
\begin{proof}
By definition of $\BL$, it is clear that
$\bl_1,\bl_2,\ldots,\bl_{w_4}$ are linearly independent.
Furthermore, from the definition of $\BM$:
\begin{equation*}
  \BM=\begin{bmatrix} \BK \\ \BL \end{bmatrix},
\end{equation*}
we find that $\bm_1,\bm_2,\ldots,\bm_{w_4}$ are
also linearly independent.
\end{proof}

We need the following lemma to prove Lemma \ref{lem5.6}:
\begin{lem}\label{lem5.5}
Suppose that
\begin{equation*}
  w_4< j_0\leq w_2 \mathrm{\ \ \ and\ \ \ }
  \bm_{j_0}\in \Col(\bm_1,\bm_2,\ldots,\bm_{j_0-1}).
\end{equation*}
Then $\bm_{w_2+1-j_0}+\bm_{j_0}$ can be expressed as
\begin{equation*}
  \bm_{w_2+1-j_0}+\bm_{j_0}=\sum_{j=w_2+2-j_0}^{j_0-1}a_j\bm_j
\end{equation*}
for some $a_j\in\ZZ/2\ZZ \ \ (j=w_2+2-j_0,\ldots,j_0-1)$.
\end{lem}
\begin{proof}
Since
\begin{equation*}
  \bm_{j_0}\in \Col(\bm_1,\bm_2,\ldots,\bm_{j_0-1}),
\end{equation*}
there are $a_j\in\ZZ/2\ZZ \ \ (j=1,2,\ldots,j_0-1)$ such that
\begin{equation*}
  \bm_{j_0}=\sum_{j=1}^{j_0-1}a_j\bm_j.
\end{equation*}
Then we have
\begin{equation*}
  \bl_{j_0}=\sum_{j=1}^{j_0-1}a_j\bl_j.
\end{equation*}

However $\BL$ has the form
\begin{equation*}
  \BL=\begin{matrix}
      \begin{matrix}
      \null & \overbrace{\hspace{1.7cm}}^{w_4=w_2/2}
            & \ \ \overbrace{\hspace{1.7cm}}^{w_4=w_2/2} \\
      \end{matrix} \\
      \begin{bmatrix}
      \bar{0} & \bar{1} & \bar{0} & \cdots & \cdots & \bar{0} & \bar{1}  \\
      \bar{0} & \bar{0} & \ddots & \bar{0} & \bar{0} & \varddots & \bar{0}  \\
      \bar{0} & \bar{0} & \bar{0}  & \bar{1} & \bar{1} & \bar{0} & \bar{0}  \\
      \bar{0} & \bar{0} & \cdots & \bar{0} & \bar{1} & \cdots & \bar{1}
      \end{bmatrix}
    \end{matrix}
\end{equation*}
or
\begin{equation*}
  \BL=\begin{matrix}
      \begin{matrix}
      \null & \overbrace{\hspace{1.6cm}}^{w_4=(w_2-1)/2} & \null \ \
            & \overbrace{\hspace{1.6cm}}^{w_4=(w_2-1)/2} \\
      \end{matrix} \\
      \begin{bmatrix}
      \bar{0} & \bar{1} & \bar{0} & \cdots & \bar{0} & \cdots & \bar{0}
        & \bar{1}  \\
      \bar{0} & \bar{0} & \ddots & \bar{0} & \bar{0} & \bar{0} & \varddots
        & \bar{0}  \\
      \bar{0} & \bar{0} & \bar{0} & \bar{1} & \bar{0} & \bar{1} & \bar{0}
        & \bar{0}  \\
      \bar{0} & \bar{0} & \bar{0} & \bar{0} & \bar{1} & \bar{0} & \bar{0}
      & \bar{0}
      \end{bmatrix}
    \end{matrix}
\end{equation*}
according as $w\equiv 0\pmod 4$ or $w\equiv 2\pmod 4$.

Thus we see
\begin{equation*}
  a_j=\bar{0} \mathrm{ \ \ \ for\ \ \ }j=1,2,\cdots,w_2-j_0
\end{equation*}
and
\begin{equation*}
  a_j=\bar{1} \mathrm{ \ \ \ for\ \ \ }j=w_2+1-j_0.
\end{equation*}
Hence we have
\begin{equation*}
  \bm_{w_2+1-j_0}+\bm_{j_0}=\sum_{j=w_2+2-j_0}^{j_0-1}a_j\bm_j.
\end{equation*}
\end{proof}

Now we are ready to prove the following lemma:
\begin{lem}\label{lem5.6}
Let
\begin{equation*}
  w=12k+2a\ \ (0<k,\ 0\leq a\leq 5),
\end{equation*}
and set
\begin{equation*}
  j_{0}=
  \begin{cases}
      4k   & \mathrm{if\ \ \ } a=0 \\
      4k+1 & \mathrm{if\ \ \ } a=1,2 \\
      4k+2 & \mathrm{if\ \ \ } a=3 \\
      4k+3 & \mathrm{if\ \ \ } a=4,5
  \end{cases} \\
      \mathrm{\ \ \ when \ \ \ }
      w<2\cdot 2^{n+2},
\end{equation*}
\begin{equation*}
  j_{0}=
  \begin{cases}
      4k   & \mathrm{if\ \ \ } a=0,1,2 \\
      4k+2 & \mathrm{if\ \ \ } a=3,4,5
  \end{cases} \\
      \mathrm{\ \ \ when \ \ \ }
      2\cdot 2^{n+2}\leq w.
\end{equation*}
Then, for $j_1$ with $w_4<j_1\leq j_0$,
we have
\begin{equation*}
  \bm_{j_1}\not\in \Col(\bm_1,\bm_2,\ldots,\bm_{j_1-1}).
\end{equation*}
\end{lem}
\begin{proof}
First we will prove the case that $j_1=j_0$
by reduction to absurdity.
So we suppose that
\begin{equation*}
  \bm_{j_0}\in \Col(\bm_1,\bm_2,\ldots,\bm_{j_0-1}).
\end{equation*}
Then, from Lemma \ref{lem5.5}, we have
\begin{equation*}
  \bm_{w_2+1-j_0}+\bm_{j_0}=\sum_{j=w_2+2-j_0}^{j_0-1}a_j\bm_j
\end{equation*}
for $a_j\in\ZZ/2\ZZ \ \, (j=w_2+2-j_0,\ldots,j_0-1)$.
Since
\begin{equation*}
  \BM=\begin{bmatrix} \BK \\ \BL \end{bmatrix},
\end{equation*}
we also have
\begin{equation*}
  \bk_{w_2+1-j_0}+\bk_{j_0}=\sum_{j=w_2+2-j_0}^{j_0-1}a_j\bk_j.
\end{equation*}
Furthermore, since $\BH$ is obtained from $\BK$ by row operations,
we have
\begin{equation}\label{eqn5.13}
  \bh_{w_2+1-j_0}+\bh_{j_0}=\sum_{j=w_2+2-j_0}^{j_0-1}a_j\bh_j.
\end{equation}

We compute that
\begin{align*}
  w_2+1-j_0
  &=
  \begin{cases}
      2k+1 & \mathrm{if\ \ \ } a=0,1 \\
      2k+2 & \mathrm{if\ \ \ } a=2,3,4 \\
      2k+3 & \mathrm{if\ \ \ } a=5
  \end{cases}
      \mathrm{\ \ \ when \ \ \ }
      w<2\cdot 2^{n+2},
\end{align*}
and that
\begin{align*}
  w_2+1-j_0
  &=
  \begin{cases}
      2k+1 & \mathrm{if\ \ \ } a=0 \\
      2k+2 & \mathrm{if\ \ \ } a=1,3 \\
      2k+3 & \mathrm{if\ \ \ } a=2,4 \\
      2k+4 & \mathrm{if\ \ \ } a=5
  \end{cases}
      \mathrm{\ \ \ when \ \ \ }
      2\cdot 2^{n+2}\leq w.
\end{align*}

Now we suppose that $w<2\cdot 2^{n+2}$,
and we take $i=\beta_1(w)$.
Then, from the identities \eqref{eqn5.7} and \eqref{eqn5.8},
we have
\begin{equation*}
  \bh_{j}(i)=\bar{0}\ \ (j=w_2+2-j_0,\ldots,j_0)
  \mathrm{\ \ \ and\ \ \ }
  \bh_{w_2+1-j_0}(i)=\bar{1}.
\end{equation*}
This contradicts \eqref{eqn5.13}.

Next we suppose that $2\cdot 2^{n+2}\leq w$,
and we take $i=\beta_2(w)$.
Then, from the identities \eqref{eqn5.9} and \eqref{eqn5.10},
we have that
\begin{equation*}
  \bh_{j}(i)=\bar{0}\ \ (j=w_2+1-j_0,\ldots,j_0-1)
  \mathrm{\ \ \ and\ \ \ }
  \bh_{j_0}(i)=\bar{1}.
\end{equation*}
This also contradicts \eqref{eqn5.13}.

Thus we have proved
\begin{equation*}
  \bm_{j_0}\not\in \Col(\bm_1,\bm_2,\ldots,\bm_{j_0-1}).
\end{equation*}

To prove
\begin{equation*}
  \bm_{j_1}\not\in \Col(\bm_1,\bm_2,\ldots,\bm_{j_1-1})
\end{equation*}
for $j_1$ with $w_4<j_1<j_0$,
we take
\begin{equation*}
  i=\beta_1(w)+2(j_0-j_1)
  \mathrm{\ \ \ if\ \ \ }
  w<2\cdot 2^{n+2}
\end{equation*}
and
\begin{equation*}
  i=\beta_2(w)-(j_0-j_1)
  \mathrm{\ \ \ if\ \ \ }
  2\cdot 2^{n+2}\leq w.
\end{equation*}
Then we apply Lemma \ref{lem5.1}.
The subsequent argument is similar to that of the case $j_1=j_0$.
The detail will be left to the reader.
\end{proof}

Combining Lemmas \ref{lem5.4} and \ref{lem5.6}, we have
\begin{lem}\label{lem5.7}
Let $j_{0}$ be as in Lemma \ref{lem5.6}.
If $\ 1<j_1\leq j_0$, then it holds that
\begin{equation*}
  \bm_{j_1}\not\in \Col(\bm_1,\bm_2,\ldots,\bm_{j_1-1}).
\end{equation*}
\end{lem}

We need the following two Lemmas \ref{lem5.8} and \ref{lem5.9}
to prove Lemma \ref{lem5.10}.
\begin{lem}\label{lem5.8}
If \ $1<j_1\leq w_2$ and $j_1$ is odd, then
\begin{equation*}
  \bm_{j_1}\not\in \Col(\bm_1,\ldots,\bm_{j_1-1}).
\end{equation*}
\end{lem}
\begin{proof}
We know
\begin{equation*}
  \bk_j(j_1-1)=\bar{0}\ \ (j=1,\ldots,j_1-1)
\end{equation*}
and
\begin{equation*}
  \bk_{j_1}(j_1-1)=\binom{j_1}{j_1-1}\mod 2=\bar{j_1}=\bar{1}.
\end{equation*}

This implies
\begin{equation*}
  \bk_{j_1}\not\in \Col(\bk_1,\ldots,\bk_{j_1-1}).
\end{equation*}
Thus we have
\begin{equation*}
  \bm_{j_1}\not\in \Col(\bm_1,\ldots,\bm_{j_1-1}).
\end{equation*}
\end{proof}

\begin{lem}\label{lem5.9}
If $w_2\equiv 0\pmod 2$, then it holds that
\begin{equation*}
  \bm_{w_2}\not\in \Col(\bm_1,\ldots,\bm_{w_2-1}).
\end{equation*}
\end{lem}
\begin{proof}
We know the second row of $\BM$ has the form
\begin{equation}\label{eqn5.14}
  (\bar{0},\bar{0},\bar{0},\bar{1},\bar{0},\bar{1},\bar{0},
  \ldots,\bar{1},\bar{0},\bar{1},\bar{0})
\end{equation}
since
\begin{equation*}
  m_{1j}
  =
  \begin{cases}
      \binom{j}{1} \mod 2 & \mathrm{if\ \ \ } 1<j \\
      \bar{0}                    & \mathrm{if\ \ \ } j\leq 1
  \end{cases}
  =
  \begin{cases}
      \bar{1} & \mathrm{if\ \ \ } 1<j,\ j\mathrm{\ \ odd} \\
      \bar{0} & \mathrm{otherwise}.
  \end{cases}
\end{equation*}

We add row vectors of $\BM$ to the vector \eqref{eqn5.14}
to obtain a vector
\begin{equation}\label{eqn5.16}
  (\bar{0},\bar{0},\bar{0},\ldots,\bar{0},\bar{0},
    \bar{1},\bar{1},\ldots,\bar{1},\bar{1},\bar{0}).
\end{equation}
Finally we add the last row vector of $\BM$
\begin{equation*}
  (\bar{0},\bar{0},\bar{0},\ldots,\bar{0},\bar{0},
    \bar{1},\bar{1},\ldots,\bar{1},\bar{1},\bar{1})
\end{equation*}
to obtain a vector
\begin{equation}\label{eqn5.17}
  (\bar{0},\bar{0},\bar{0},\ldots,\bar{0},\bar{0},\bar{0},\bar{0},
    \ldots,\bar{0},\bar{0},\bar{1}).
\end{equation}
This shows that the vector \eqref{eqn5.17} is
a linear combination of the row vectors of $\BM$.
Hence we know that the last column vector $\bm_{w_2}$
of $\BM$ is linearly independent of
$\bm_{j}\ (j=1,2,\ldots,w_2-1)$.
This completes the proof.
\end{proof}

Now we define sets of column vectors in $\BM$
as follows:
\begin{defn}\label{defn5.3}
Let
\begin{equation*}
  w=12k+2a\ \ (0<k,\ 0\leq a\leq 5).
\end{equation*}
\begin{enumerate}
\item
If $a=0$, we set
\begin{equation*}
  S:=
  \{\bm_1,\bm_2,\ldots,\bm_{4k-1},\bm_{4k}\}\cup
  \{\bm_{4k+1},\bm_{4k+1+2},\ldots,\bm_{4k+1+2(k-1)}\}\cup
  \{\bm_{6k}\};
\end{equation*}
\item
If $a=1$, we set
\begin{equation*}
  S:=
  \{\bm_1,\bm_2,\ldots,\bm_{4k-1},\bm_{4k}\}\cup
  \{\bm_{4k+1},\bm_{4k+1+2},\ldots,\bm_{4k+1+2k}\};
\end{equation*}
\item
If $a=2$, we set
\begin{equation*}
  S:=
  \{\bm_1,\bm_2,\ldots,\bm_{4k-1},\bm_{4k}\}\cup
  \{\bm_{4k+1},\bm_{4k+1+2},\ldots,\bm_{4k+1+2k}\}\cup
  \{\bm_{6k+2}\};
\end{equation*}
\item
If $a=3$, we set
\begin{equation*}
  S:=
  \{\bm_1,\bm_2,\ldots,\bm_{4k+1},\bm_{4k+2}\}\cup
  \{\bm_{4k+3},\bm_{4k+3+2},\ldots,\bm_{4k+3+2k}\};
\end{equation*}
\item
If $a=4$, we set
\begin{equation*}
  S:=
  \{\bm_1,\bm_2,\ldots,\bm_{4k+1},\bm_{4k+2}\}\cup
  \{\bm_{4k+3},\bm_{4k+3+2},\ldots,\bm_{4k+3+2k}\}\cup
  \{\bm_{6k+4}\};
\end{equation*}
\item
If $a=5$, we set
\begin{equation*}
  S:=
  \{\bm_1,\bm_2,\ldots,\bm_{4k+1},\bm_{4k+2}\}\cup
  \{\bm_{4k+3},\bm_{4k+3+2},\ldots,\bm_{4k+3+2(k+1)}\}.
\end{equation*}
\end{enumerate}
\end{defn}

Now, concluding the discussion in this section, we have
\begin{lem}\label{lem5.10}
The vectors in $S$ are linearly independent.
\end{lem}
\begin{proof}
Suppose that $j>1$ and $\bm_j\in S$.
Then,
if $w_2\equiv 0\pmod 2$, by Lemmas \ref{lem5.7}, \ref{lem5.8}
and \ref{lem5.9},
and if $w_2\equiv 1\pmod 2$, by Lemmas \ref{lem5.7} and
\ref{lem5.8},
we know
\begin{equation*}
  \bm_{j}\not\in \Col(\bm_1,\bm_2,\ldots,\bm_{j-1}).
\end{equation*}
This implies that $S$ is a set of linearly independent vectors.
\end{proof}

\section{Proofs of Theorems \ref{thm1.2} and \ref{thm1.3}}
\label{sect6}

In this section, we give proofs of Theorems \ref{thm1.2}
and \ref{thm1.3}. A straightforward computation gives:
\begin{lem}\label{lem6.1}
The cardinality of $S$, say \#S, satisfies
\begin{equation*}
  w_2-\#S=d_w.
\end{equation*}
\end{lem}

Now we are ready to prove Theorems \ref{thm1.2} and \ref{thm1.3}.
\begin{proof}[Proof of Theorem \ref{thm1.2}]
One can easily check that Theorem \ref{thm1.2} holds for $w_2$
with $w_2<12$.
Therefore we assume that $12\leq w_2$
as in Assumption \ref{assumpt4.7}.

Noting that $\#S=w_2-d_w$ by Lemma \ref{lem6.1},
we can express $S$ as follows:
\begin{enumerate}
\item
If $w\equiv 0\pmod 4$, then
\begin{equation*}
  S=\{\bm_j\ |\ 1\leq j\leq w_2,\  j\ne w_2-2,w_2-4,\ldots,w_2-2d_w\};
\end{equation*}
\item
If $w\equiv 2\pmod 4$, then
\begin{equation*}
  S=\{\bm_j\ |\ 1\leq j\leq w_2,\  j\ne w_2-1,w_2-1-2,\ldots,w_2-1-2(d_w-1)\}.
\end{equation*}
\end{enumerate}

Now let $n_w$ be the integer defined by
\begin{equation*}
  n_w:=w_2-d_w,
  \text{\ \ that is,\ \ \ }
  n_w=\#S,
\end{equation*}
and let
\begin{equation*}
  j_a\ \ (a=0,1,\cdots,n_w-1)
\end{equation*}
be integers such that
\begin{equation*}
  1\leq j_0<j_1<\cdots<j_{n_w-1}\leq w_2
  \text{\ \ \ and\ \ \ }
  \bm_{j_a}\in S.
\end{equation*}
Furthermore, let
\begin{equation*}
  \hat{j}_b\ \ (b=0,1,\cdots,d_w-1)
\end{equation*}
be integers such that
\begin{equation*}
  1\leq \hat{j}_0<\hat{j}_1<\cdots<\hat{j}_{d_w-1}\leq w_2
  \text{\ \ \ and\ \ \ }
  \bm_{\hat{j}_b}\not\in S.
\end{equation*}
Clearly, we have
\begin{equation*}
  S=\{\bm_{j_a}\ |\ a=0,1,\cdots,n_w-1\}
\end{equation*}
and
\begin{equation}\label{eqn6.1}
\begin{split}
 \{ \hat{j}_0,& \hat{j}_1, \ldots, \hat{j}_{d_w-1}\} \\
 &=\begin{cases}
 \{ w_2-2,w_2-4,\ldots,w_2-2d_w\}
         & \mathrm{if\ \ \ } w\equiv 0 \pmod 4 \\
 \{ w_2-1,w_2-1-2,\ldots,w_2-1-2(d_w-1)\}
         & \mathrm{if\ \ \ } w\equiv 2 \pmod 4.
    \end{cases}
\end{split}
\end{equation}

Next we consider matrices $\BD_1$ and $\BM_1$ defined by
\begin{align*}
  \BD_1&:=[\bd_{j_0},\bd_{j_1},\ldots,\bd_{j_{n_w-1}}], \\
  \BM_1&:=[\bm_{j_0},\bm_{j_1},\ldots,\bm_{j_{n_w-1}}].
\end{align*}
Note that $\BM_1$ is the reduction modulo 2 of $\BD_1$.

Since the column vectors
$\bm_{j_0},\bm_{j_1},\ldots,\bm_{j_{n_w-1}}$
are linearly independent, we can choose $n_w$ rows,
say rows $i_0,i_1,\ldots,i_{n_w-1}$ of $\BM_1$, so that
the matrix $\BM_2$ defined by
\begin{equation*}
  \BM_2:=[m_{{i_c}{j_a}}]\ \ (c,a=0,1,\ldots,n_w-1)
\end{equation*}
has non-zero determinant.
Then the matrix $\BD_2$ defined by
\begin{equation*}
  \BD_2:=[d_{{i_c}{j_a}}]\ \ (c,a=0,1,\ldots,n_w-1)
\end{equation*}
also has non-zero determinant.
This is because $\BM_2$ is the reduction modulo 2 of $\BD_2$.

Furthermore, we set
\begin{equation*}
  \BD_3:=[d_{{i_c}{\hat{j}_b}}]\ \ (c=0,1,\ldots,n_w-1;\ b=0,1,\ldots,d_w-1).
\end{equation*}
Then, from Lemma \ref{lem2.3}, we know that
\begin{equation*}
    \BD_2
    \begin{bmatrix}
      t_{w,j_0} \\ t_{w,j_1} \\ \vdots \\ t_{w,j_{n_w-1}}
    \end{bmatrix}
  +
    \BD_3
    \begin{bmatrix}
      t_{w,\hat{j}_0} \\ t_{w,\hat{j}_1} \\ \vdots \\ t_{w,\hat{j}_{d_w-1}}
    \end{bmatrix}
  =\zr,
\end{equation*}
from which it follows that
\begin{equation*}
    \begin{bmatrix}
      t_{w,j_0} \\ t_{w,j_1} \\ \vdots \\ t_{w,j_{n_w-1}}
    \end{bmatrix}
  =
    -\BD_2^{-1}\BD_3
    \begin{bmatrix}
      t_{w,\hat{j}_0} \\ t_{w,\hat{j}_1} \\ \vdots \\ t_{w,\hat{j}_{d_w-1}}
    \end{bmatrix}.
\end{equation*}
This implies that  each of $\ t_{w,j_0},t_{w,j_1},\ldots,t_{w,j_{n_w-1}}\ $
can be expressed as linear combinations of
$\ t_{w,\hat{j}_0},t_{w,\hat{j}_1},\ldots,t_{w,\hat{j}_{d_w-1}}\ $
over $\QQ$.
On the other hand, by Theorem \ref{thm1.1}, we know that
$\ r_{w,1},r_{w,3},\ldots,r_{w,w-1}\ $ span $S_{w+2}^*$.
Thus
$\ t_{w,1},t_{w,2},\ldots,t_{w,w_2}\ $ also span $S_{w+2}^*$.
These imply that
\begin{equation}\label{eqn6.2}
  t_{w,\hat{j}_0},t_{w,\hat{j}_1},\ldots,t_{w,\hat{j}_{d_w-1}}
  \text{\ \ \ span\ \ \ } S_{w+2}^*.
\end{equation}

From the fact \eqref{eqn6.2} together with the identity \eqref{eqn6.1},
and the identities
\begin{equation*}
  t_{w,j}=s_{w,2j-1}
  \text{\ \ and\ \ }
  s_{w,n}=(-1)^n\binom{w}{n}r_{w,w-n},
\end{equation*}
it follows that,
\begin{align*}
 \{r_{w,w-(2j-1)}\ &|\ j=w_2-2,w_2-4,\ldots,w_2-2d_w\} \\
 &=\{r_{w,4i+1}\ |\ i=1,2,\ldots,d_w\}
\end{align*}
span $S_{w+2}^*$,
when $w\equiv 0 \pmod 4$,
and
\begin{align*}
 \{r_{w,w-(2j-1)}\ &|\ j=w_2-1,w_2-3,\ldots,w_2-1-2(d_w-1)\} \\
 &=\{r_{w,4i-1}\ |\ i=1,2,\ldots,d_w\}
\end{align*}
span $S_{w+2}^*$,
when $w\equiv 2 \pmod 4$.

Since
\begin{equation*}
  \dim S_{w+2}^*=d_w,
\end{equation*}
we know that
\begin{equation*}
 \{r_{w,4i\pm1}\ |\ i=1,2,\ldots,d_w\}
\end{equation*}
is a basis of $S_{w+2}^*$.
This completes the proof of Theorem \ref{thm1.2}.
\end{proof}

\begin{proof}[Proof of Theorem \ref{thm1.3}]
The cusp form $R_{w,m}$ is characterized
by the formula:
\begin{equation*}
  r_{w,m}(f)=(R_{w,m},f)\ \mathrm{\ \ \ for\ any \ \ } f\in S_{w+2}
\end{equation*}
with the Petersson inner product $(\ ,\ )$.
Then, from Theorem \ref{thm1.2}, it follows that
\begin{equation*}
  R_{w,{4i\pm 1}}\ (i=1,2,\ldots,d_w)
\end{equation*}
is a basis for $S_{w+2}$.
This completes the proof.
\end{proof}

\section{Relationship between $R_{w,n}$, $E_{w,n}$ and $S_{w,n}$}
\label{sect7}

The Sections from \ref{sect7} to \ref{sect13}
are devoted to the study of Hecke operators.
We also present proofs for Theorems \ref{thm1.6}, \ref{thm1.7},
\ref{thm1.8} and \ref{thm1.9}.

First we show that, for $n$ odd,
 $\alpha_{w+2}^+$ maps $c_{w,n}R_{w,n}$ to $E_{w,n}$,
and that $\beta_{w}^+$ maps $E_{w,n}$ to $S_{w,n}$.
From this we will prove
Theorems \ref{thm1.6} and \ref{thm1.7}.

It is known (\cite[Theorem $1'$]{KZ1}) that
 $R_{w,n}\ (0<n< w;\ n \text{\ odd})$
corresponds to
 $S_{w,n}$
by the Eichler-Shimura monomorphism
 $\beta_{w}^+\alpha_{w+2}^+:S_{w+2}\to\mathcal{U}_w^+$.
More precisely we know that
\begin{equation}\label{eqn7.1}
  \beta_{w}^+\alpha_{w+2}^+(c_{w,n}R_{w,{n}})=S_{w,{n}}.
\end{equation}
It was also shown (\cite[Theorem 5.2]{F4}) that
 $E_{w,n}\ (0<n< w;\ n \text{\ odd})$
corresponds to
 $S_{w,n}$
by the map
 $\beta_{w}^+:\mathcal{E}_w^+\to\mathcal{U}_w^+$.  Namely,
\begin{equation}\label{eqn7.2}
  \beta_{w}^+(E_{w,{n}})=S_{w,{n}}.
\end{equation}

From \eqref{eqn7.1} and \eqref{eqn7.2}, we obtain:
\begin{lem}\label{lem7.1}
Let $n$ be an odd integer with  $0<n<w$.
Then it holds that
\begin{equation}\label{eqn7.3}
  \alpha_{w+2}^+(c_{w,n}R_{w,{n}})=E_{w,{n}}.
\end{equation}
\end{lem}
\begin{proof}
Since $\ker \beta_{w}^+$ is spanned by $G_w$,
from \eqref{eqn7.1} and \eqref{eqn7.2},
we know that
\begin{equation*}
  \alpha_{w+2}^+(c_{w,n}R_{w,{n}})-E_{w,{n}}=cG_w
\end{equation*}
for some constant $c$.

On the other hand, we have
\begin{align*}
  \alpha_{w+2}^+(c_{w,n}R_{w,n})(1,0)
    &=\int_{0}^{\infty}c_{w,n}R_{w,n}(z)z^wdz
      \text{\ \ \ (by \eqref{eqn1.1})} \\
    &=\beta_{w}^+\alpha_{w+2}^+(c_{w,n}R_{w,n})(1,0)
      \text{\ \ \ (by \eqref{eqn1.2})} \\
    &=S_{w,n}(1,0)
      \text{\ \ \ (by \eqref{eqn7.1})} \\
    &=-\frac{B_{n+1}}{n+1}
      -\frac{B_{\tn+1}}{\tn+1}
      +\frac{w+2}{B_{w+2}}\frac{B_{n+1}}{n+1}\frac{B_{\tn+1}}{\tn+1} \\
    &=E_{w,n}(1,0).
\end{align*}
Hence we know the constant $c$ is zero,
from which we arrive at \eqref{eqn7.3}.
\end{proof}

Now we are ready to prove Theorems \ref{thm1.6} and \ref{thm1.7}.
\begin{proof}[Proofs of Theorems \ref{thm1.6} and \ref{thm1.7}]
The Eichler-Shimura theorem (\cite[pp.\ 199--200]{KZ1})
asserts that the map
\begin{equation*}
  \beta_{w}^+\alpha_{w+2}^+:S_{w+2}\to
  \mathcal{E}_w^+\to\mathcal{U}_w^+
\end{equation*}
is an monomorphism and that the image
$\beta_{w}^+\alpha_{w+2}^+(S_{w+2})$
and $h^w-k^w$ span $\mathcal{U}_w^+$.
Since $R_{w,{4i\pm 1}}\ (i=1,2,\ldots,d_w)$ is a basis of
$S_{w+2}$ and
$\beta_{w}^+\alpha_{w+2}^+(c_{w,n}R_{w,{n}})=S_{w,{n}}$,
we know that $S_{w,{4i\pm 1}}\ (i=1,2,\ldots,d_w)$ and
$h^w-k^w$ is a basis of $\mathcal{U}_w^+$.
This completes the proof of Theorem \ref{thm1.6}.

Next, by Theorem \ref{thm1.4}, we see that the map
\begin{equation*}
  \alpha_{w+2}^+:S_{w+2}\to
  \mathcal{E}_w^+
\end{equation*}
is an monomorphism and that the image
$\alpha_{w+2}^+(S_{w+2})$, $F_w$ and $G_w$
span $\mathcal{E}_w^+$. Furthermore, from Lemma \ref{lem7.1},
we know that
\begin{equation*}
  \alpha_{w+2}^+(c_{w,4i\pm 1}R_{w,{4i\pm 1}})
  =E_{w,{4i\pm 1}}.
\end{equation*}
Hence we see that $E_{w,{4i\pm 1}}\ (i=1,2,\ldots,d_w)$,
$F_w$ and $G_w$ is a basis of $\mathcal{U}_w^+$.
This completes the proof of Theorem \ref{thm1.7}.
\end{proof}

Next we show that, for $n$ even,
 $\alpha_{w+2}^-$ maps $c_{w,n}R_{w,n}$ to $E_{w,n}$.
For even $n$ with $0<n<w$, it was shown
(\cite[Theorem $1'$]{KZ1}, \cite[Theorem 5.2]{F4}) that
\begin{equation}\label{eqn7.4}
  \beta_{w}^-\alpha_{w+2}^-(c_{w,n}R_{w,{n}})=S_{w,{n}}
  \text{\ \ \ and\ \ \ \ }
  \beta_{w}^-(E_{w,{n}})=S_{w,{n}}.
\end{equation}

From the identities \eqref{eqn7.4} and the fact that
$\beta_{w}^-$ is an isomorphism, we obtain:
\begin{lem}\label{lem7.2}
Let $n$ be an even integer with  $0<n<w$.
Then it holds that
\begin{equation}\label{eqn7.6}
  \alpha_{w+2}^-(c_{w,n}R_{w,{n}})=E_{w,{n}}.
\end{equation}
\end{lem}

\section{The action of the Hecke operator on $R_{w,n}$}
\label{sect8}

In Definition \ref{defn1.2}, we introduced $R_{w,n}^m$
as a generalization of $R_{w,n}$, i.e. $R_{w,n}^1=R_{w,n}$.
Naturally $R_{w,n}^m$ inherits most of the properties from $R_{w,n}$.
Moreover we can understand the action of the Hecke operator
on $R_{w,n}$ in terms of $R_{w,n}^m$.
Indeed, we will show that $T_mR_{w,n}=R_{w,n}^m$.

For this purpose, we need the following notation.
Let $\gamma$ be a matrix such that
\begin{equation*}
       \gamma=\begin{bmatrix}a&b\\c&d\end{bmatrix}
       \text{\ \ with $ad-bc>0$\ \ },
\end{equation*}
and let
\begin{equation*}
       f:\ \{z\in\CC\ |\ \Im z>0\}\ \to\ \CC
\end{equation*}
be a map. Then, as usual, we define $f|\gamma$ by
\begin{equation*}
       (f|\gamma)(z):=(cz+d)^{-w-2}f(\frac{az+b}{cz+d}).
\end{equation*}
Let $M_m$ be the complete set of right coset
representatives of $H_m$ modulo $\Gamma$, defined as follows:
\begin{equation*}
  M_m:=\left\{
       \begin{bmatrix}a&b\\0&d\end{bmatrix}
       \ |\ \ ad=m,\ \ a>0,\ \ b \mod d
       \right\}.
\end{equation*}
Recall that the action of the Hecke operation $T_m$ on $f$ is given by
\begin{equation*}
       T_mf(z):=
         m^{w+1}\sum_{\alpha \in M_m}(f|\alpha)(z).
\end{equation*}

Now we are ready to prove the following lemma:
\begin{lem}\label{lem8.1}
\begin{equation*}
  T_mR_{w,n}=R_{w,n}^m.
\end{equation*}
\end{lem}
\begin{proof}
We have
\begin{align*}
  T_mR_{w,n}(z)
    &=c_{w,n}^{-1}T_m\sum_{\substack{\left[\begin{smallmatrix}a&b\\c&d\end{smallmatrix}\right]
                      \in \Gamma}}
           \frac{1}{(az+b)^{n+1}(cz+d)^{\tn+1}} \\
    &=c_{w,n}^{-1}T_m\sum_{\gamma \in \Gamma}
           z^{-n-1}|\gamma \\
    &=m^{w+1}c_{w,n}^{-1}\sum_{\alpha \in M_m}\sum_{\gamma \in \Gamma}
           z^{-n-1}|\gamma|\alpha \\
    &=m^{w+1}c_{w,n}^{-1}\sum_{\gamma' \in H_m}
           z^{-n-1}|\gamma' \\
    &=m^{w+1}c_{w,n}^{-1}\sum_{\substack
             {\left[\begin{smallmatrix}a'&b'\\c'&d'\end{smallmatrix}\right]
                      \in H_m}}
           \frac{1}{(a'z+b')^{n+1}(c'z+d')^{\tn+1}} \\
    &=R_{w,n}^{m}(z).
\end{align*}
Here we used the identity
\begin{equation*}
  H_m=\bigcup_{\alpha \in M_m}\Gamma \alpha.
\end{equation*}

This completes the proof.
\end{proof}

\section{Properties of $E_{w,n}^m$}
\label{sect9}

In this section, we study $E_{w,n}^m$ which generalizes $E_{w,n}$.
We will show that $E_{w,n}^m$ is a well-defined Dedekind
symbols of weight $w$.
Furthermore we will show that $E_{w,n}^m$ is even or odd
depending on $n$ is odd or even.

First we will show that the sum in Definition \ref{defn1.2} (3)
is finite by the following two lemmas
(the proofs are left to the reader).
\begin{lem}\label{lem9.1}
Let $h$ be a positive integer.
Suppose that
$\left[\begin{smallmatrix}a&b\\c&d\end{smallmatrix}\right]
          \in H_m$
satisfies $ac\ne 0$ and
\begin{equation}\label{eqn9.1}
  (\frac{k}{h}+\frac{b}{a})(\frac{k}{h}+\frac{d}{c})<0.
\end{equation}
Then
\begin{equation*}
  mh\geq |a|,\ \ mh\geq |c|.
\end{equation*}
\end{lem}

\begin{lem}\label{lem9.2}
Let $h$ be a positive integer.
We consider two conditions for
$\left[\begin{smallmatrix}a&b\\c&d\end{smallmatrix}\right]
\in H_m$:
\begin{equation}\label{eqn9.2}
  ac\ne 0;
\end{equation}
\begin{equation}\label{eqn9.3}
  (\frac{k}{h}+\frac{b}{a})(\frac{k}{h}+\frac{d}{c})<0.
\end{equation}
Then there are only finite many
$\left[\begin{smallmatrix}a&b\\c&d\end{smallmatrix}\right]
\in H_m$
which satisfy \eqref{eqn9.2} and \eqref{eqn9.3}.
\end{lem}

By Lemmas \ref{lem9.1} and \ref{lem9.2}, we know that
$E_{w,n}^m$ is well-defined since the sum
in Definition \ref{defn1.2} (3)
is finite.

Next we will show that $E_{w,n}^m$ is a Dedekind symbol
in the following lemma:
\begin{lem}\label{lem9.3}
$E_{w,n}^m$ is an even $($resp. odd$)$ Dedekind symbol
of weight $w$ for $n$ odd $($resp. even$)$.
\end{lem}
\begin{proof}
The argument used in \cite[Lemma 3.1 and Theorem 1.1]{F3}
for $D_{w,n}$ is also valid to establish that $E_{w,n}^m$
is an odd or even Dedekind symbol according as $n$ is even or odd.
The reader should refer to \cite{F3} for the detail.

It is clear that $E_{w,n}^m$ satisfies
\begin{equation*}
  E_{w,n}^m(ch,ck)=c^wE_{w,n}^m(h,k)
  \text{\ \ \ for any\ \ } c\in\ZZ^+.
\end{equation*}
Hence the weight of $E_{w,n}^m$ is $w$.
This completes the proof.
\end{proof}

\section{The action of the Hecke operator on $E_{w,n}$}
\label{sect10}

In this section we will study how $T_m$ acts on $E_{w,n}$.
We will prove the following proposition.
\begin{prop}\label{prop10.1}
Let $m$ and $n$ be positive integers such that $0<n<w$.
Then it hold that
\begin{equation*}
    T_mE_{w,n}=E_{w,n}^m.
\end{equation*}
\end{prop}
\begin{proof}
Let $E_{w,n,0}^m$ and $E_{w,n,0}$ be defined by
\begin{align*}
  E_{w,n,0}^m:&=\frac{1}{2}
    \sum_{\substack{\left[\begin{smallmatrix}a&b\\c&d\end{smallmatrix}\right]
                      \in H_m \\ ac\ne 0 \\
      (k/h+b/a)(k/h+d/c)<0}}
    \sgn(\frac{k}{h}+\frac{b}{a})(ak+bh)^{\tn}(ck+dh)^n, \\
  E_{w,n,0}:&=\frac{1}{2}
    \sum_{\substack{\left[\begin{smallmatrix}a&b\\c&d\end{smallmatrix}\right]
                      \in \Gamma \\ ac\ne 0 \\
      (k/h+b/a)(k/h+d/c)<0}}
    \sgn(\frac{k}{h}+\frac{b}{a})(ak+bh)^{\tn}(ck+dh)^n. \\
\end{align*}
First we will prove that $T_mE_{w,n,0}=E_{w,n,0}^m$.

Noting that
\begin{equation*}
  H_m=\bigcup_{\alpha \in M_m}\Gamma \alpha
    =\bigcup_{\substack{ad=m \\ d>0 \\ b\, (\md \, d)}}\Gamma
    \begin{bmatrix}a&b\\0&d\end{bmatrix},
\end{equation*}
we have
{\allowdisplaybreaks
\begin{align*}
  T_m&E_{w,n,0}(h,k) \\
    &=\sum_{\substack{ad=m \\ d>0}}\sum_{b\,(\md\, d)}
    E_{w,n,0}(dh,ak+bh) \\
    &=\frac{1}{2}\sum_{\substack{ad=m \\ d>0}}\sum_{b\,(\md\, d)}
    \sum_{\substack
         {\left[\begin{smallmatrix}a'&b'\\c'&d'\end{smallmatrix}\right]
                      \in \Gamma \\ a'c'\ne 0 \\
         \{(ak+bh)/dh+b'/a'\}\{(ak+bh)/dh+d'/c'\}<0}}
    \sgn(\frac{ak+bh}{dh}+\frac{b'}{a'}) \\
    &\cccsp\ccsp\times\{a'(ak+bh)+b'dh\}^{\tn}\{c'(ak+bh)+d'dh\}^n \\
    &=\frac{1}{2}\sum_{\substack{ad=m \\ d>0}}\sum_{b\,(\md\, d)}
    \sum_{\substack
         {\left[\begin{smallmatrix}a'&b'\\c'&d'\end{smallmatrix}\right]
                      \in \Gamma \\ a'c'\ne 0 \\
         \{k/h+(a'b+b'd)/a'a\}\{k/h+(c'b+d'd)/c'a\}<0}}
    \sgn(\frac{k}{h}+\frac{a'b+b'd}{a'a}) \\
    &\cccsp\ccsp\times\{a'ak+(a'b+b'd)h\}^{\tn}\{c'ak+(c'b+d'd)h\}^n \\
    &=\frac{1}{2}
    \sum_{\substack
         {\left[\begin{smallmatrix}a''&b''\\c''&d''\end{smallmatrix}\right]
                      \in H_m \\ a''c''\ne 0 \\
         (k/h+b''/a'')(k/h+d''/c'')<0}}
    \sgn(\frac{k}{h}+\frac{b''}{a''})
    (a''k+b''h)^{\tn}(c''k+d''h)^n \\
    &=E_{w,n,0}^m(h,k)
\end{align*}
}
where we set
\begin{equation*}
\left\{
\begin{aligned}
    a''&=a'a \\
    b''&=a'b+b'd \\
    c''&=c'a \\
    d''&=c'b+d'd,
\end{aligned}
\right.
\text{\ \ \ \ \ \ namely,\ \ \ }
   \begin{bmatrix}a''&b''\\c''&d''\end{bmatrix}
   =\begin{bmatrix}a'&b'\\c'&d'\end{bmatrix}
   \begin{bmatrix}a&b\\0&d\end{bmatrix}.
\end{equation*}

Next, applying the formula
\begin{equation*}
  \sum_{b=0}^{c-1}\bar{B}_{n+1}(x+\frac{b}{c})
    =c^{-n}\bar{B}_{n+1}(cx)\ \ \ (c\in\ZZ^+),
\end{equation*}
we obtain the other terms of $T_mE_{w,n}(h,k)$ as follows:
{\allowdisplaybreaks
\begin{align*}
  T_m\frac{\bar{B}_{n+1}(\frac{k}{h})h^w}{n+1}
    &=\sum_{\substack{ad=m \\ d>0}}\sum_{b\,(\md\, d)}
      \frac{\bar{B}_{n+1}(\frac{ak+bh}{dh})(dh)^w}{n+1} \\
    &=\sum_{\substack{ad=m \\ d>0}}
      \frac{d^{-n}\bar{B}_{n+1}(d\frac{ak}{dh})(dh)^w}{n+1} \\
    &=\sum_{\substack{ad=m \\ d>0}}
      d^{\tn}\frac{\bar{B}_{n+1}(\frac{ak}{h})h^w}{n+1}, \\
  T_m\frac{\bar{B}_{\tn+1}(\frac{k}{h})h^w}{\tn+1}
    &=\sum_{\substack{ad=m \\ d>0}}
      d^{n}\frac{\bar{B}_{\tn+1}(\frac{ak}{h})h^w}{\tn+1}
\end{align*}
}
and
\begin{align*}
  T_mh^w
    &=\sum_{\substack{ad=m \\ d>0}}\sum_{b\,(\md\, d)}(dh)^w
     =\sum_{\substack{ad=m \\ d>0}}d^{w+1}h^w
     =\sigma_{w+1}(m)h^w.
\end{align*}
From these identities we obtain $T_mE_{w,n}=E_{w,n}^m$
completing the proof.
\end{proof}

\section{Reciprocity law for Dedekind symbols $E_{w,n}^m$}
\label{sect11}

In this section we will prove
the following reciprocity law for Dedekind symbols $E_{w,n}^m$:

\begin{prop}\label{prop11.1}
Let $m$ and $n$ be positive integers such that $0<n<w$.
Then it holds that
\begin{equation*}
    E_{w,n}^m(h,k)-E_{w,n}^m(k,-h)=S_{w,n}^m(h,k)
\end{equation*}
for any $(h,k)\in \ZZ^+\times\ZZ^+$.
\end{prop}

We need a few lemmas to prove Proposition \ref{prop11.1}.
We consider the sum in Definition \ref{defn1.2}
and express this as the sum of three series:
\begin{equation*}
  \begin{split}
  \frac{1}{2}
    \sum_{\substack{\left[\begin{smallmatrix}a&b\\c&d\end{smallmatrix}\right]
                      \in H_m \\ ac\ne 0 \\
      (k/h+b/a)(k/h+d/c)<0}}
    \sgn(\frac{k}{h}+\frac{b}{a})&(ak+bh)^{\tn}(ck+dh)^n \\
  &=E_{w,n,1}^m(h,k)+E_{w,n,2}^m(h,k)+E_{w,n,3}^m(h,k)
  \end{split}
\end{equation*}
where
\begin{equation*}
  E_{w,n,1}^m(h,k):=\frac{1}{2}
    \sum_{\substack{\left[\begin{smallmatrix}a&b\\c&d\end{smallmatrix}\right]
                      \in H_m \\ abcd>0 \\
      (k/h+b/a)(k/h+d/c)<0}}
    \sgn(\frac{k}{h}+\frac{b}{a})(ak+bh)^{\tn}(ck+dh)^n,
\end{equation*}
\begin{equation*}
  E_{w,n,2}^m(h,k):=\frac{1}{2}
    \sum_{\substack{\left[\begin{smallmatrix}a&b\\c&d\end{smallmatrix}\right]
                      \in H_m \\ abcd<0 \\
      (k/h+b/a)(k/h+d/c)<0}}
    \sgn(\frac{k}{h}+\frac{b}{a})(ak+bh)^{\tn}(ck+dh)^n
\end{equation*}
and
\begin{equation*}
  E_{w,n,3}^m(h,k):=\frac{1}{2}
    \sum_{\substack{\left[\begin{smallmatrix}a&b\\c&d\end{smallmatrix}\right]
                      \in H_m \\ ac\ne 0,\ bd=0 \\
      (k/h+b/a)(k/h+d/c)<0}}
    \sgn(\frac{k}{h}+\frac{b}{a})(ak+bh)^{\tn}(ck+dh)^n.
\end{equation*}

First we show a reciprocal property of $E_{w,n,1}^m$.
\begin{lem}\label{lem11.2}
Suppose that $(h,k)\in \ZZ^+\times\ZZ^+$.
Then
\begin{equation*}
    E_{w,n,1}^m(h,k)= E_{w,n,1}^m(k,-h).
\end{equation*}
\end{lem}
\begin{proof}
From the definition of $E_{w,n,1}^m(h,k)$, we have
\begin{equation*}
  E_{w,n,1}^m(h,k)=\frac{1}{2}
    \sum_{\substack{\left[\begin{smallmatrix}a&b\\c&d\end{smallmatrix}\right]
                      \in H_m \\ abcd>0 \\
      (k/h+b/a)(k/h+d/c)<0}}
    \sgn(\frac{k}{h}+\frac{b}{a})(ak+bh)^{\tn}(ck+dh)^n.
\end{equation*}
Replacing
$\left[\begin{smallmatrix}a&b\\c&d\end{smallmatrix}\right]$
with
$\left[\begin{smallmatrix}b&-a\\d&-c\end{smallmatrix}\right]$,
we can transform this into the following formula:

\begin{equation}\label{eqn11.1}
  \begin{split}
  \frac{1}{2}
    &\sum_{\substack{\left[\begin{smallmatrix}a&b\\c&d\end{smallmatrix}\right]
                      \in H_m \\ abcd>0 \\
      (k/h+b/a)(k/h+d/c)<0}}
    \sgn(\frac{k}{h}+\frac{b}{a})(ak+bh)^{\tn}(ck+dh)^n \\
  &\ccsp=\frac{1}{2}
    \sum_{\substack{\left[\begin{smallmatrix}a&b\\c&d\end{smallmatrix}\right]
                      \in H_m \\ abcd>0 \\
      (k/h-a/b)(k/h-c/d)<0}}
    \sgn(\frac{k}{h}-\frac{a}{b})(bk-ah)^{\tn}(dk-ch)^n.
  \end{split}
\end{equation}
On the other hand, we have
\begin{equation}\label{eqn11.2}
  E_{w,n,1}^m(k,-h)=
    \frac{1}{2}
    \sum_{\substack{\left[\begin{smallmatrix}a&b\\c&d\end{smallmatrix}\right]
                      \in H_m \\ abcd>0 \\
      (-h/k+b/a)(-h/k+d/c)<0}}
    \sgn(-\frac{h}{k}+\frac{b}{a})(-ah+bk)^{\tn}(-ch+dk)^n.
\end{equation}

Now we compare the right hand sides of
\eqref{eqn11.1} and \eqref{eqn11.2}.
First we see that the condition $abcd>0$
implies that $\sgn(ac)=\sgn(bd)$ and $\sgn(ab)=\sgn(cd)$.
Next we see that
  $$h^2bd(k/h-a/b)(k/h-c/d)=(-ah+bk)(-ch+dk)
  =k^2ac(-h/k+b/a)(-h/k+d/c).$$
From these we know that two conditions
$(k/h-a/b)(k/h-c/d)<0$
and
$(-h/k+b/a)(-h/k+d/c)<0$
are equivalent.
Furthermore, under the condition that
$(k/h-a/b)(k/h-c/d)<0$, we get
$0<k/h<a/b$ or $0<k/h<c/d$,
and then $ab>0$ or $cd>0$.
This together with $\sgn(ab)=\sgn(cd)$ implies
$\sgn(ab)>0$, and then $\sgn(a)=\sgn(b)$.
Hence we have
  $$\sgn(\frac{k}{h}-\frac{a}{b})
  =\sgn(hb)\sgn(bk-ah)
  =\sgn(ka)\sgn(bk-ah)
  =\sgn(-\frac{h}{k}+\frac{b}{a}).$$
In conclusion we see that the right hand sides
of \eqref{eqn11.1} and \eqref{eqn11.2}
are equal.
This completes the proof.
\end{proof}

Next we show a reciprocal property of $E_{w,n,2}^m$.
\begin{lem}\label{lem11.3}
Suppose that $(h,k)\in \ZZ^+\times\ZZ^+$.
Then
\begin{equation*}
    E_{w,n,2}^m(h,k)-E_{w,n,2}^m(k,-h)=\frac{1}{2}
    \sum_{\substack{\left[\begin{smallmatrix}a&b\\c&d\end{smallmatrix}\right]
                      \in H_m \\ abcd<0 }}
    \sgn(ab)(ak+bh)^{\tn}(ck+dh)^n.
\end{equation*}
\end{lem}
\begin{proof}
From the definition of $E_{w,n,2}^m(h,k)$, we have
{\allowdisplaybreaks
\begin{align*}
  E&_{w,n,2}^m(h,k) \\
    &=\frac{1}{2}
    \sum_{\substack{\left[\begin{smallmatrix}a&b\\c&d\end{smallmatrix}\right]
                      \in H_m \\ abcd<0 \\ (k/h+b/a)(k/h+d/c)<0}}
    \sgn(\frac{k}{h}+\frac{b}{a})(ak+bh)^{\tn}(ck+dh)^n \\
    &=\frac{1}{2}
    \sum_{\substack{\left[\begin{smallmatrix}a&b\\c&d\end{smallmatrix}\right]
                      \in H_m \\ abcd<0 \\
                      -d/c>k/h>-b/a}}
    (ak+bh)^{\tn}(ck+dh)^n
    -\frac{1}{2}
    \sum_{\substack{\left[\begin{smallmatrix}a&b\\c&d\end{smallmatrix}\right]
                      \in H_m \\ abcd<0 \\
                      -d/c<k/h<-b/a}}
    (ak+bh)^{\tn}(ck+dh)^n \\
    &=\frac{1}{2}
    \sum_{\substack{\left[\begin{smallmatrix}a&b\\c&d\end{smallmatrix}\right]
                      \in H_m \\ abcd<0 \\
                      -d/c>k/h}}
    (ak+bh)^{\tn}(ck+dh)^n
    -\frac{1}{2}
    \sum_{\substack{\left[\begin{smallmatrix}a&b\\c&d\end{smallmatrix}\right]
                      \in H_m \\ abcd<0 \\
                      k/h<-b/a}}
    (ak+bh)^{\tn}(ck+dh)^n \\
    &\cccsp \text{(noting $\sgn(-d/c)=-\sgn(-b/a)$)}.
\end{align*}
}
On the other hand,
replacing
$\left[\begin{smallmatrix}a&b\\c&d\end{smallmatrix}\right]$
with
$\left[\begin{smallmatrix}b&-a\\d&-c\end{smallmatrix}\right]$,
we have
\begin{align*}
  -&E_{w,n,2}^m(k,-h) \\
    &=-\frac{1}{2}
    \sum_{\substack{\left[\begin{smallmatrix}a&b\\c&d\end{smallmatrix}\right]
                      \in H_m \\ abcd<0 \\ (-h/k+b/a)(-h/k+d/c)<0}}
    \sgn(\frac{-h}{k}+\frac{b}{a})(-ah+bk)^{\tn}(-ch+dk)^n \\
    &=-\frac{1}{2}
    \sum_{\substack{\left[\begin{smallmatrix}a&b\\c&d\end{smallmatrix}\right]
                      \in H_m \\ abcd<0 \\ (-h/k-a/b)(-h/k-c/d)<0}}
    \sgn(\frac{-h}{k}+\frac{-a}{b})(-bh-ak)^{\tn}(-dh-ck)^n \\
    &\cccsp \text{(replacing
        $\left[\begin{smallmatrix}a&b\\c&d\end{smallmatrix}\right]$
        with
        $\left[\begin{smallmatrix}b&-a\\d&-c\end{smallmatrix}\right]$)} \\
    &=-\frac{1}{2}
    \sum_{\substack{\left[\begin{smallmatrix}a&b\\c&d\end{smallmatrix}\right]
                      \in H_m \\ abcd<0 \\ (-h/k-a/b)(-h/k-c/d)<0}}
    \sgn(\frac{-h}{k}+\frac{-a}{b})(ak+bh)^{\tn}(ck+dh)^n \\
    &=\frac{1}{2}
    \sum_{\substack{\left[\begin{smallmatrix}a&b\\c&d\end{smallmatrix}\right]
                      \in H_m \\ abcd<0 \\
                      -c/d>h/k>-a/b}}
    (ak+bh)^{\tn}(ck+dh)^n
    -\frac{1}{2}
    \sum_{\substack{\left[\begin{smallmatrix}a&b\\c&d\end{smallmatrix}\right]
                      \in H_m \\ abcd<0 \\
                      -c/d<h/k<-a/b}}
    (ak+bh)^{\tn}(ck+dh)^n \\
    &=\frac{1}{2}
    \sum_{\substack{\left[\begin{smallmatrix}a&b\\c&d\end{smallmatrix}\right]
                      \in H_m \\ abcd<0 \\
                      -c/d>h/k}}
    (ak+bh)^{\tn}(ck+dh)^n
    -\frac{1}{2}
    \sum_{\substack{\left[\begin{smallmatrix}a&b\\c&d\end{smallmatrix}\right]
                      \in H_m \\ abcd<0 \\
                      h/k<-a/b}}
    (ak+bh)^{\tn}(ck+dh)^n \\
    &=\frac{1}{2}
    \sum_{\substack{\left[\begin{smallmatrix}a&b\\c&d\end{smallmatrix}\right]
                      \in H_m \\ abcd<0 \\
                      0<-d/c<k/h}}
    (ak+bh)^{\tn}(ck+dh)^n
    -\frac{1}{2}
    \sum_{\substack{\left[\begin{smallmatrix}a&b\\c&d\end{smallmatrix}\right]
                      \in H_m \\ abcd<0 \\
                      k/h>-b/a>0}}
    (ak+bh)^{\tn}(ck+dh)^n. \\
\end{align*}
Hence we have
{\allowdisplaybreaks
\begin{align*}
  E_{w,n,2}^m(h,k)-&E_{w,n,2}^m(k,-h) \\
    &=\frac{1}{2}
    \sum_{\substack{\left[\begin{smallmatrix}a&b\\c&d\end{smallmatrix}\right]
                      \in H_m \\ abcd<0 \\
                      -d/c>k/h}}
    -\frac{1}{2}
    \sum_{\substack{\left[\begin{smallmatrix}a&b\\c&d\end{smallmatrix}\right]
                      \in H_m \\ abcd<0 \\
                      k/h<-b/a}}
    +\frac{1}{2}
    \sum_{\substack{\left[\begin{smallmatrix}a&b\\c&d\end{smallmatrix}\right]
                      \in H_m \\ abcd<0 \\
                      0<-d/c<k/h}}
    -\frac{1}{2}
    \sum_{\substack{\left[\begin{smallmatrix}a&b\\c&d\end{smallmatrix}\right]
                      \in H_m \\ abcd<0 \\
                      k/h>-b/a>0}} \\
    &=\frac{1}{2}
    \sum_{\substack{\left[\begin{smallmatrix}a&b\\c&d\end{smallmatrix}\right]
                      \in H_m \\ abcd<0 \\
                      \sgn(-d/c)>0}}
    -\frac{1}{2}
    \sum_{\substack{\left[\begin{smallmatrix}a&b\\c&d\end{smallmatrix}\right]
                      \in H_m \\ abcd<0 \\
                      \sgn(-b/a)>0}} \\
    &=\frac{1}{2}
    \sum_{\substack{\left[\begin{smallmatrix}a&b\\c&d\end{smallmatrix}\right]
                      \in H_m \\ abcd<0 \\
                      \sgn(ab)>0}}
    -\frac{1}{2}
    \sum_{\substack{\left[\begin{smallmatrix}a&b\\c&d\end{smallmatrix}\right]
                      \in H_m \\ abcd<0 \\
                      \sgn(ab)<0}} \\
    &=\frac{1}{2}
    \sum_{\substack{\left[\begin{smallmatrix}a&b\\c&d\end{smallmatrix}\right]
                      \in H_m \\ abcd<0}}
                      \sgn(ab)(ak+bh)^{\tn}(ck+dh)^n.
\end{align*}
}
This completes the proof.
\end{proof}

Finally we look into the term $E_{w,n,3}^m(h,k)$.
We prove that $E_{w,n,3}^m(h,k)$ is expressed in terms
of Bernoulli polynomials and Bernoulli functions.
\begin{lem}\label{lem11.4}
Suppose that
$h$ and $k$ are positive integers.
Then
\begin{equation}\label{eqn11.3}
  \begin{split}
    E_{w,n,3}^m(h,k)&=
      \sum_{\substack{ad=m \\ a>0}}
    d^{\tn}
      \left\{
        \frac{B_{n+1}(\frac{ah}{k})}{n+1}
        -\frac{\bar{B}_{n+1}(\frac{ah}{k})}{n+1}
      \right\}k^{w} \\
    &\csp-(-1)^{n}
      \sum_{\substack{ad=m \\ a>0}}
    d^{n}
      \left\{
        \frac{B_{\tn+1}(\frac{ah}{k})}{\tn+1}
        -\frac{\bar{B}_{\tn+1}(\frac{ah}{k})}{\tn+1}
      \right\}k^{w}.
  \end{split}
\end{equation}
\end{lem}
\begin{proof}
From the definition of $E_{w,n,3}^m$, we have
\begin{equation}\label{eqn11.4}
\begin{split}
  E_{w,n,3}^m(h,k)=&\frac{1}{2}
    \sum_{\substack{\left[\begin{smallmatrix}a&b\\c&d\end{smallmatrix}\right]
                      \in H_m \\ ac\ne 0,\ b=0 \\
      (k/h+b/a)(k/h+d/c)<0}}
    \sgn(\frac{k}{h}+\frac{b}{a})(ak+bh)^{\tn}(ck+dh)^n \\
    &+\frac{1}{2}
    \sum_{\substack{\left[\begin{smallmatrix}a&b\\c&d\end{smallmatrix}\right]
                      \in H_m \\ ac\ne 0,\ d=0 \\
      (k/h+b/a)(k/h+d/c)<0}}
    \sgn(\frac{k}{h}+\frac{b}{a})(ak+bh)^{\tn}(ck+dh)^n.
  \end{split}
\end{equation}
We first consider the first sum in the right hand side of \eqref{eqn11.4}.
We will apply the following formula
\begin{equation}\label{eqn11.5}
  \frac{B_{m+1}(x+1)}{m+1}-\frac{B_{m+1}(x)}{m+1}=x^m\ \ (m\geq 0)
\end{equation}
to transform the sum into the following formula:
{\allowdisplaybreaks
\begin{align*}
  \frac{1}{2}
    &\sum_{\substack{\left[\begin{smallmatrix}a&b\\c&d\end{smallmatrix}\right]
                      \in H_m \\ ac\ne 0,\ b=0 \\
      (k/h+b/a)(k/h+d/c)<0}}
    \sgn(\frac{k}{h}+\frac{b}{a})(ak+bh)^{\tn}(ck+dh)^n \\
    &\csp=
      \sum_{\substack{\left[\begin{smallmatrix}a&b\\c&d\end{smallmatrix}\right]
                      \in H_m \\ ac\ne 0,\ b=0,\ a>0 \\
      (k/h+b/a)(k/h+d/c)<0}}
    \sgn(\frac{k}{h}+\frac{b}{a})(ak+bh)^{\tn}(ck+dh)^n \\
    &\csp=
      \sum_{\substack{ad=m,\ a>0 \\ c\ne 0 \\ (k/h+d/c)<0}}
    (ak)^{\tn}(ck+dh)^n \\
    &\csp=
      \sum_{\substack{ad=m \\ a>0}}
      \sum_{-dh/k<c<0}
    k^wa^{\tn}\left(c+\frac{dh}{k}\right)^n \\
    &\csp=
      \sum_{\substack{ad=m \\ a>0}}
      \sum_{dh/k>c>0}
    k^wa^{\tn}\left(\frac{dh}{k}-c\right)^n \\
    &\csp=
      \sum_{\substack{ad=m \\ a>0}}
      \sum_{c=1}^{\lfl\frac{dh}{k}\rfl}
    k^wa^{\tn}\left(\frac{dh}{k}-c\right)^n \\
    &\csp=
      \sum_{\substack{ad=m \\ a>0}}
    k^{w}a^{\tn}
      \left\{
        \frac{B_{n+1}(\frac{dh}{k})}{n+1}
        -\frac{B_{n+1}(\frac{dh}{k}-\lfl\frac{dh}{k}\rfl)}{n+1}
      \right\} \text{\ (applying \eqref{eqn11.5} repeatedly)} \\
    &\csp=
      \sum_{\substack{ad=m \\ a>0}}
    a^{\tn}
      \left\{
        \frac{B_{n+1}(\frac{dh}{k})}{n+1}
        -\frac{\bar{B}_{n+1}(\frac{dh}{k})}{n+1}
      \right\}k^{w} \\
    &\csp=
      \sum_{\substack{ad=m \\ a>0}}
    d^{\tn}
      \left\{
        \frac{B_{n+1}(\frac{ah}{k})}{n+1}
        -\frac{\bar{B}_{n+1}(\frac{ah}{k})}{n+1}
      \right\}k^{w}.
\end{align*}
}
Similarly the second sum can be transformed into the following formula:
  \begin{align*}
  \frac{1}{2}
    \sum_{\substack{\left[\begin{smallmatrix}a&b\\c&d\end{smallmatrix}\right]
                      \in H_m \\ ac\ne 0,\ d=0 \\
      (k/h+b/a)(k/h+d/c)<0}}
    &\sgn(\frac{k}{h}+\frac{b}{a})(ak+bh)^{\tn}(ck+dh)^n  \\
    &=
    -(-1)^{n}
      \sum_{\substack{ad=m \\ a>0}}
    d^{n}
      \left\{
        \frac{B_{\tn+1}(\frac{ah}{k})}{\tn+1}
        -\frac{\bar{B}_{\tn+1}(\frac{ah}{k})}{\tn+1}
      \right\}k^w.
  \end{align*}
Combining these we obtain \eqref{eqn11.3}.
This completes the proof.
\end{proof}

Similar computation yields the following result:
\begin{lem}\label{lem11.5}
Suppose that
$h$ and $k$ are positive integers.
Then
\begin{align*}
    E_{w,n,3}^m(k,-h)&=
      \sum_{\substack{ad=m \\ a>0}}
    d^{n}
      \left\{
        \frac{B_{\tn+1}(\frac{ak}{h})}{\tn+1}
        -\frac{\bar{B}_{\tn+1}(\frac{ak}{h})}{\tn+1}
      \right\}h^{w} \\
    &\csp+(-1)^{\tn+1}
    \sum_{\substack{ad=m \\ a>0}}
    d^{\tn}
      \left\{
        \frac{B_{n+1}(\frac{ak}{h})}{n+1}
        -\frac{\bar{B}_{n+1}(\frac{ak}{h})}{n+1}
      \right\}h^{w}.
\end{align*}
\end{lem}

Now we are ready to prove Proposition \ref{prop11.1}.
\begin{proof}
[\it{Proof of Proposition \ref{prop11.1}}]
From Lemmas \ref{lem11.2}, \ref{lem11.3}, \ref{lem11.4} and \ref{lem11.5},
we have
{\allowdisplaybreaks
\begin{align*}
  E_{w,n}^m&(h,k)-E_{w,n}^m(k,-h) \\
  =&E_{w,n,1}^m(h,k)-E_{w,n,1}^m(k,-h)+E_{w,n,2}^m(h,k)-E_{w,n,2}^m(k,-h) \\
    &+E_{w,n,3}^m(h,k)-E_{w,n,3}^m(k,-h) \\
      &+
      \sum_{\substack{ad=m \\ a>0}}
      \left\{(-1)^n
        d^{\tn}\frac{\bar{B}_{n+1}(\frac{ak}{h})h^w}{n+1}
        -d^n\frac{\bar{B}_{\tn+1}(\frac{ak}{h})h^w}{\tn+1}
        \right\} \\
      &+
      \sum_{\substack{ad=m \\ a>0}}
      \left\{
        d^{\tn}\frac{\bar{B}_{n+1}(\frac{ah}{k})k^w}{n+1}
        -(-1)^nd^n\frac{\bar{B}_{\tn+1}(\frac{ah}{k})k^w}{\tn+1}
        \right\} \\
      &+
      \begin{cases}
      \sigma_{w+1}(m)\frac{w+2}{B_{w+2}}
      \frac{B_{n+1}}{n+1}\frac{B_{\tn+1}}{\tn+1}(h^w-k^w)
         & \mathrm{if\ \ \ } n\equiv 1 \pmod {2} \\
      0  & \mathrm{if\ \ \ } n\equiv 0 \pmod {2}
      \end{cases} \\
  =&\frac{1}{2}
    \sum_{\substack{\left[\begin{smallmatrix}a&b\\c&d\end{smallmatrix}\right]
                      \in H_m \\ abcd<0 }}
    \sgn(ab)(ak+bh)^{\tn}(ck+dh)^n \\
    &+
      \sum_{\substack{ad=m \\ a>0}}
    d^{\tn}
      \left\{
        \frac{B_{n+1}(\frac{ah}{k})}{n+1}
        -\frac{\bar{B}_{n+1}(\frac{ah}{k})}{n+1}
      \right\}k^{w} \\
    &\csp-(-1)^{n}
      \sum_{\substack{ad=m \\ a>0}}
    d^{n}
      \left\{
        \frac{B_{\tn+1}(\frac{ah}{k})}{\tn+1}
        -\frac{\bar{B}_{\tn+1}(\frac{ah}{k})}{\tn+1}
      \right\}k^{w} \\
    &-
      \sum_{\substack{ad=m \\ a>0}}
    d^{n}
      \left\{
        \frac{B_{\tn+1}(\frac{ak}{h})}{\tn+1}
        -\frac{\bar{B}_{\tn+1}(\frac{ak}{h})}{\tn+1}
      \right\}h^{w} \\
    &\csp-(-1)^{\tn+1}
    \sum_{\substack{ad=m \\ a>0}}
    d^{\tn}
      \left\{
        \frac{B_{n+1}(\frac{ak}{h})}{n+1}
        -\frac{\bar{B}_{n+1}(\frac{ak}{h})}{n+1}
      \right\}h^{w} \\
      &+
      \sum_{\substack{ad=m \\ a>0}}
      \left\{(-1)^n
        d^{\tn}\frac{\bar{B}_{n+1}(\frac{ak}{h})h^w}{n+1}
        -d^n\frac{\bar{B}_{\tn+1}(\frac{ak}{h})h^w}{\tn+1}
        \right\} \\
      &+
      \sum_{\substack{ad=m \\ a>0}}
      \left\{
        d^{\tn}\frac{\bar{B}_{n+1}(\frac{ah}{k})k^w}{n+1}
        -(-1)^nd^n\frac{\bar{B}_{\tn+1}(\frac{ah}{k})k^w}{\tn+1}
        \right\} \\
      &+
      \begin{cases}
      \sigma_{w+1}(m)\frac{w+2}{B_{w+2}}
      \frac{B_{n+1}}{n+1}\frac{B_{\tn+1}}{\tn+1}(h^w-k^w)
         & \mathrm{if\ \ \ } n\equiv 1 \pmod {2} \\
      0  & \mathrm{if\ \ \ } n\equiv 0 \pmod {2}
      \end{cases} \\
  =&\frac{1}{2}
      \sum_{\substack{\left[\begin{smallmatrix}a&b\\c&d\end{smallmatrix}\right]
                      \in H_m \\
                      abcd<0}}
      \sgn(ab)(ak+bh)^{\tn}(ck+dh)^{n} \\
      &+\sum_{\substack{ad=m \\ a>0}}
        \left\{
        (-1)^nd^\tn\frac{B_{n+1}(\frac{ak}{h})h^w}{n+1}
        +d^\tn\frac{B_{n+1}(\frac{ah}{k})k^w}{n+1}
        \right. \\
      &\cccsp\ccsp\ \
        \left.
        -d^n\frac{B_{\tn+1}(\frac{ak}{h})h^w}{\tn+1}
        -(-1)^nd^n\frac{B_{\tn+1}(\frac{ah}{k})k^w}{\tn+1}
        \right\} \\
      &+
      \begin{cases}
      \sigma_{w+1}(m)\frac{w+2}{B_{w+2}}
      \frac{B_{n+1}}{n+1}\frac{B_{\tn+1}}{\tn+1}(h^w-k^w)
         & \mathrm{if\ \ \ } n\equiv 1 \pmod {2} \\
      0  & \mathrm{if\ \ \ } n\equiv 0 \pmod {2}
      \end{cases} \\
  =&S_{w,n}^m(h,k).
\end{align*}
}
Thus we obtain the required reciprocity laws completing the proof of
Proposition \ref{prop11.1}.
\end{proof}

\section{The action of the Hecke operator on $S_{w,n}$}
\label{sect12}

Now we arrive at the following proposition
which will be the key of obtaining our explicit formulas
for Hecke operators:
\begin{prop}\label{prop12.1}
Let $m$ be a positive integer, and let
$n$ be an integer such that $0<n<w$.
Then it hold that
\begin{equation*}
  T_m(S_{w,{n}})=S_{w,{n}}^m.
\end{equation*}
\end{prop}
\begin{proof}
We notice that the Hecke operators on the spaces of cusp forms,
the spaces of Dedekind symbols and the spaces of
period polynomials are all compatible
(the diagrams \eqref{eqn1.3} and \eqref{eqn1.4}).

First we prove the case that $n$ is odd.
We have
{\allowdisplaybreaks
\begin{align*}
  T_m(S_{w,n})
  &=T_m(\beta_{w}^+\alpha_{w+2}^+ (c_{w,n}R_{w,n}))
    \ \ \text{(by  \eqref{eqn7.1})} \\
  &=\beta_{w}^+\alpha_{w+2}^+T_m (c_{w,n}R_{w,n})
    \ \ \text{(from the diagram \eqref{eqn1.3})} \\
  &=\beta_{w}^+T_m\alpha_{w+2}^+ (c_{w,n}R_{w,n})
    \ \ \text{(from the diagram \eqref{eqn1.4})} \\
  &=\beta_{w}^+T_m(E_{w,n})
    \ \ \text{(by Lemma \ref{lem7.1})} \\
  &=\beta_{w}^+(E_{w,n}^m)
    \ \ \text{(by Proposition \ref{prop10.1})} \\
  &=S_{w,n}^m
    \ \ \text{(by Proposition \ref{prop11.1})}.
\end{align*}
}
Substituting
Lemma \ref{lem7.2}, $\alpha_{w+2}^-$ and $\beta_{w}^-$
for
Lemma \ref{lem7.1}, $\alpha_{w+2}^+$ and $\beta_{w}^+$,
we can prove the case that $n$ is even.
This completes the proof.
\end{proof}

We are also ready to prove Theorem \ref{thm1.8}.
\begin{proof}[Proof of Theorem \ref{thm1.8}]
It is obvious that Theorem \ref{thm1.8} follows
from Lemma \ref{lem8.1}, Proposition \ref{prop10.1}
and Proposition \ref{prop12.1}.
\end{proof}

\section{Explicit formulas for Hecke operators}
\label{sect13}

In this section, we give a proof to Theorem \ref{thm1.9}.
\begin{proof}[Proof of Theorem \ref{thm1.9}]

We suppose that
\begin{equation*}
  \BT_{m}=[\tau_{ij}]\ \ (i,j=1,2,\ldots,d_w)
\end{equation*}
is the matrix representing the Hecke operator
\ $T_m:S_{w+2}\to S_{w+2}$
with respect to the basis
\begin{equation*}
  c_{w,4i\pm 1}R_{w,{4i\pm 1}}\ \ (i=1,2,\ldots,d_w).
\end{equation*}
Namely, we suppose that
\begin{equation}\label{eqn13.1}
  T_m(c_{w,4j\pm 1}R_{w,{4j\pm 1}})
    =\sum_{i=1}^{d_w}\tau_{ij}c_{w,4i\pm 1}R_{w,{4i\pm 1}}
  \ \ (j=1,2,\ldots,d_w).
\end{equation}

Then we have
\begin{align*}
    S_{w,{4j\pm 1}}^m
    &=T_m(S_{w,{4j\pm 1}})
    \ \ \text{(by Proposition \ref{prop12.1})} \\
    &=T_m\beta_{w}^+\alpha_{w+2}^+ (c_{w,4j\pm 1}R_{w,{4j\pm 1}})
    \ \ \text{(by \eqref{eqn7.1})} \\
    &=\beta_{w}^+\alpha_{w+2}^+ T_m(c_{w,4j\pm 1}R_{w,{4j\pm 1}})
    \ \ \text{(from the diagram \eqref{eqn1.3})} \\
    &=\sum_{i=1}^{d_w}\tau_{ij}
     \beta_{w}^+\alpha_{w+2}^+ (c_{w,4i\pm 1}R_{w,{4i\pm 1}})
    \ \ \text{(by \eqref{eqn13.1})} \\
    &=\sum_{i=1}^{d_w}\tau_{ij}S_{w,{4i\pm 1}}
    \ \ \text{(by \eqref{eqn7.1})}.
\end{align*}
This gives
\begin{equation}\label{eqn13.2}
    S_{w,{4j\pm 1}}^m
     =\sum_{i=1}^{d_w}\tau_{ij}S_{w,{4i\pm 1}}
    \ \ (j=1,2,\ldots,d_w).
\end{equation}
Furthermore, by taking inner product of $S_{w,{4k\pm 1}}$
and each side of \eqref{eqn13.2},
we have
\begin{align*}
  \langle S_{w,{4k\pm 1}},S_{w,{4j\pm 1}}^m\rangle
  &=\sum_{i=1}^{d_w}
  \langle S_{w,{4k\pm 1}},S_{w,{4i\pm 1}}\rangle
  \tau_{ij}
  \ \ \ (j,k=1,2,\ldots,d_w).
\end{align*}
Namely, we have
\begin{align*}
   \begin{bmatrix}
   \langle S_{w,{4i\pm 1}},S_{w,{4j\pm 1}}^m\rangle
   \end{bmatrix}=
   \begin{bmatrix}
   \langle S_{w,{4i\pm 1}},S_{w,{4j\pm 1}}\rangle
   \end{bmatrix}
   \BT_{m}
   \ \ \ (i,j=1,2,\ldots,d_w).
\end{align*}
Now the linear independence of
$S_{w,{4i\pm 1}}\ (i=1,2,\ldots,d_w)$
guarantees that the matrix
\begin{align*}
   \begin{bmatrix}
   \langle S_{w,{4i\pm 1}},S_{w,{4j\pm 1}}\rangle
   \end{bmatrix}
   \ \ \ (i,j=1,2,\ldots,d_w)
\end{align*}
is non-singular.
Thus we have
\begin{equation*}
  \BT_{m}=
   \begin{bmatrix}
   \langle S_{w,{4i\pm 1}},S_{w,{4j\pm 1}}\rangle
   \end{bmatrix}^{-1}
   \begin{bmatrix}
   \langle S_{w,{4i\pm 1}},S_{w,{4j\pm 1}}^m\rangle
   \end{bmatrix}
   \ \ (i,j=1,2,\ldots,d_w)
\end{equation*}
completing the proof of (1).

Next we will prove (2).
First we expand the term
\begin{equation*}
    \frac{1}{2}
    \sum_{\substack{\left[\begin{smallmatrix}a&b\\c&d\end{smallmatrix}\right]
                      \in H_m \\ abcd<0 }}
    \sgn(ab)(ak+bh)^{\tn}(ck+dh)^n
\end{equation*}
as a polynomial in $h$ and $k$ as follows:
{\allowdisplaybreaks
\begin{align*}
    \frac{1}{2}&
    \sum_{\substack{\left[\begin{smallmatrix}a&b\\c&d\end{smallmatrix}\right]
                      \in H_m \\ abcd<0 }}
    \sgn(ab)(ak+bh)^{\tn}(ck+dh)^n \\
    &=
    \sum_{\substack{\left[\begin{smallmatrix}a&b\\c&d\end{smallmatrix}\right]
                      \in H_m \\ a>0\ abcd<0 }}
    \sgn(ab)(ak+bh)^{\tn}(ck+dh)^n  \\
    &=
    \sum_{\substack{\left[\begin{smallmatrix}a&b\\c&d\end{smallmatrix}\right]
                      \in H_m \\ a>0\ b>0\ abcd<0 }}
    (ak+bh)^{\tn}(ck+dh)^n
    -\sum_{\substack{\left[\begin{smallmatrix}a&b\\c&d\end{smallmatrix}\right]
                      \in H_m \\ a>0\ b<0\ abcd<0 }}
    (ak+bh)^{\tn}(ck+dh)^n  \\
    &=
    \sum_{\substack{ad-bc=m \\ a>0\ b>0\ abcd<0 }}
    (ak+bh)^{\tn}(ck+dh)^n
    -\sum_{\substack{ad-bc=m \\ a>0\ b<0\ abcd<0 }}
    (ak+bh)^{\tn}(ck+dh)^n  \\
    &=
    \sum_{\mu=1}^{m-1}
    \sum_{\substack{ad=\mu \\ a>0}}
    \sum_{\substack{bc=\mu-m \\ b>0}}
    (ak+bh)^{\tn}(ck+dh)^n
    -\sum_{\mu=1}^{m-1}
    \sum_{\substack{ad=\mu \\ a>0}}
    \sum_{\substack{bc=\mu-m \\ b>0}}
    (ak-bh)^{\tn}(-ck+dh)^n \\
    &=
    2\sum_{\mu=1}^{m-1}
    \sum_{\substack{ad=\mu \\ a>0}}
    \sum_{\substack{bc=\mu-m \\ b>0}}
    \sum_{\substack{0\leq k\leq \tn\\ 0\leq \ell\leq n\\ k+n-\ell \text{\ odd}}}
    \binom{\tn}{k}a^{\tn-k}k^{\tn-k}b^{k}h^{k}
    \binom{n}{\ell}c^{n-\ell} k^{n-\ell} d^{\ell}h^{\ell} \\
    &=
    2\sum_{\mu=1}^{m-1}
    \sum_{\substack{ad=\mu \\ 0<a}}
    \sum_{\substack{bc=\mu-m \\ 0<b}}
    \sum_{\substack{0\leq k\leq \tn\\ 0\leq \ell\leq n\\ k+\ell \text{\ even}}}
    \binom{\tn}{k}a^{\tn-k}k^{\tn-k}b^{k}h^{k}
    \binom{n}{\ell}(\frac{\mu-m}{b})^{n-\ell} k^{n-\ell}
      (\frac{\mu}{a})^{\ell}h^{\ell} \\
    &=
    2\sum_{\mu=1}^{m-1}
    \sum_{\substack{ad=\mu \\ 0<a}}
    \sum_{\substack{bc=\mu-m \\ 0<b}}
    \sum_{\substack{0\leq k\leq \tn\\ 0\leq \ell\leq n\\ k+\ell \text{\ even}}}
    \mu^{\ell}(\mu-m)^{n-\ell}
    \binom{\tn}{k}\binom{n}{\ell}
    a^{\tn-k-\ell}b^{k+\ell-n}h^{k+\ell}k^{\tn+n-k-\ell} \\
    &=
    2\sum_{\mu=1}^{m-1}
    \sum_{\substack{ad=\mu \\ 0<a}}
    \sum_{\substack{bc=\mu-m \\ 0<b}}
    \sum_{\substack{\nu=0 \\ \nu\text{\ even}}}^{w}
    \sum_{\lambda=\max(0,\nu-\tn)}^{\min(n,\nu)}
    \mu^{\lambda}(\mu-m)^{n-\lambda}
    \binom{\tn}{\nu-\lambda}\binom{n}{\lambda}
    a^{\tn-\nu}b^{\nu-n}h^{\nu}k^{w-\nu} \\
    &\cccsp\text{(setting $\nu=k+\ell$ and $\lambda=\ell$)} \\
    &=
    2\sum_{\mu=1}^{m-1}
    \sum_{\substack{\nu=0 \\ \nu\text{\ even}}}^{w}
    \sum_{\lambda=\max(0,\nu-\tn)}^{\min(n,\nu)}
    \mu^{\lambda}(\mu-m)^{n-\lambda}
    \binom{\tn}{\nu-\lambda}\binom{n}{\lambda}
    \sigma_{\tn-\nu}(\mu)\sigma_{\nu-n}(m-\mu)
    h^{\nu}k^{w-\nu} \\
    &=
    2\sum_{\substack{\nu=0 \\ \nu\text{\ even}}}^{w}
    \sum_{\mu=1}^{m-1}
    \sum_{\lambda=\max(0,\nu-\tn)}^{\min(n,\nu)}
    \mu^{\lambda}(\mu-m)^{n-\lambda}
    \binom{\tn}{\nu-\lambda}\binom{n}{\lambda}
    \sigma_{\tn-\nu}(\mu)\sigma_{\nu-n}(m-\mu)
    h^{\nu}k^{w-\nu}. \\
\end{align*}
}
We also calculate other terms:
{\allowdisplaybreaks
\begin{align*}
      \sum_{\substack{ad=m \\ a>0}}
        d^\tn\frac{B_{n+1}(\frac{ak}{h})h^w}{n+1}
      &=\sum_{\substack{ad=m \\ a>0}}
        \frac{d^\tn}{n+1}\sum_{\mu=0}^{n+1}\binom{n+1}{\mu}
        B_{\mu}\left(\frac{ak}{h}\right)^{n+1-\mu}h^w \\
      &=\sum_{\substack{ad=m \\ a>0}}
        \frac{m^{\tn}}{n+1}\sum_{\mu=0}^{n+1}\binom{n+1}{\mu}
        B_{\mu}a^{n+1-\mu-\tn}h^{w-n-1+\mu}k^{n+1-\mu} \\
      &=
        \frac{m^\tn}{n+1}\sum_{\mu=0}^{n+1}\binom{n+1}{\mu}
        B_{\mu}\sigma_{n+1-\mu-\tn}(m)h^{w-n-1+\mu}k^{n+1-\mu} \\
      &=
        \frac{m^\tn}{n+1}\sum_{\nu=\tn-1}^{w}\binom{n+1}{\nu-\tn+1}
        B_{\nu-\tn+1}\sigma_{n-\nu}(m)h^{\nu}k^{w-\nu}, \\
\end{align*}
}
\begin{align*}
      \sum_{\substack{ad=m \\ a>0}}
        d^\tn\frac{B_{n+1}(\frac{ah}{k})k^w}{n+1}
      &=
        \frac{m^\tn}{n+1}\sum_{\nu=0}^{n+1}\binom{n+1}{n-\nu+1}
        B_{n-\nu+1}\sigma_{\nu-\tn}(m)h^{\nu}k^{w-\nu}. \\
\end{align*}

Summing up these identities, we obtain the formula for $S_{w,n}^m$
in (2) of Theorem \ref{thm1.9}.
This completes the proof.
\end{proof}

\begin{center}
\bf{Appendix : Computing matrices representing the Hecke operators
and their characteristic polynomials}
\end{center}
\label{appendix}

Finally in this section, we will demonstrate a
Mathematica\footnote{Mathematica is a trademark of Wolfram Research,~Inc.}
program which, for given $w$
and $m$, yields a matrix representing
the Hecke operator $T_m$ on $S_{w+2}$,
and its characteristic polynomial.
The program is a straightforward application of Theorem \ref{thm1.9}
using built-in functions BernoulliB and DivisorSigma.

\vspace{3mm}
\begin{verbatim}
w=28; (* set a positive even integer *)
m=7; (* set a positive integer *)
If[Mod[w,12]==0,dw=Quotient[w+2,12]-1,dw=Quotient[w+2,12]];
If[OddQ[w]||(dw<=0)||(m<=0),
  Print["w odd or dw<=0 or m<=0"];Exit[]];
swn=Array[a,{2,dw,w+1}];
For[j=1,j<=2,j++,
  If[j<=1,em=1,em=m];
  For[i=1,i<=dw,i++,
    If[Mod[w,4]==0,n=4*i+1,n=4*i-1]; tn=w-n;
    For[nu=1,nu<=w+1,nu++,
      swn[[j,i,nu]]=0;
      lmin=Min[n,nu-1]; lmax=Max[0,nu-1-tn];
      If[OddQ[nu],
        For[mu=1,mu<=em-1,mu++,
          For[lambda=lmax,lambda<=lmin,lambda++,
            swn[[j,i,nu]]=
              swn[[j,i,nu]]+2*(mu^lambda)*((mu-em)^(n-lambda))*
              Binomial[tn,nu-1-lambda]*Binomial[n,lambda]*
              DivisorSigma[tn-nu+1,mu]*DivisorSigma[nu-1-n,em-mu];
          ];
        ];
      ];
    ];
    For[nu=tn,nu<=w+1,nu++,
      swn[[j,i,nu]]=swn[[j,i,nu]]-(em^tn)*Binomial[n+1,nu-tn]*
        BernoulliB[nu-tn]*DivisorSigma[n+1-nu,em]/(n+1);
    ];
    For[nu=n,nu<=w+1,nu++,
      swn[[j,i,nu]]=swn[[j,i,nu]]-(em^n)*Binomial[tn+1,nu-n]*
        BernoulliB[nu-n]*DivisorSigma[tn+1-nu,em]/(tn+1);
    ];
    For[nu=1,nu<=n+2,nu++,
      swn[[j,i,nu]]=swn[[j,i,nu]]+(em^tn)*Binomial[n+1,n-nu+2]*
        BernoulliB[n-nu+2]*DivisorSigma[nu-1-tn,em]/(n+1);
    ];
    For[nu=1,nu<=tn+2,nu++,
      swn[[j,i,nu]]=swn[[j,i,nu]]+(em^n)*Binomial[tn+1,tn-nu+2]*
        BernoulliB[tn-nu+2]*DivisorSigma[nu-1-n,em]/(tn+1);
    ];
    swn[[j,i,w+1]]=
      swn[[j,i,w+1]]+DivisorSigma[w+1,em]*((w+2)*BernoulliB[n+1]*
      BernoulliB[tn+1])/(BernoulliB[w+2]*(n+1)*(tn+1));
    swn[[j,i,1]]=
      swn[[j,i,1]]-DivisorSigma[w+1,em]*((w+2)*BernoulliB[n+1]*
      BernoulliB[tn+1])/(BernoulliB[w+2]*(n+1)*(tn+1));
  ];
];
s1=Array[b,{dw,dw}];
s2=Array[c,{dw,dw}];
For[i=1,i<=dw,i++,
  For[j=1,j<=dw,j++,s1[[i,j]]=swn[[1,i]].swn[[1,j]];
    s2[[i,j]]=swn[[1,i]].swn[[2,j]]]];
twm=Inverse[s1].s2;
Print["representation matrix=",MatrixForm[twm]];
Print["characteristic polynomial=",Det[x*IdentityMatrix[dw]-twm]];
\end{verbatim}

\vspace{3mm}
For example, the output in the case that $w=28$ and $m=7$ is as follows:

\vspace{3mm}
  representation matrix=
  $\begin{pmatrix}
    \frac{-597428921326303528}{6439} &
    \frac{-4321468293778944}{6439} \\
    \frac{79904984173167605760}{6439} &
    \frac{577981127961754328}{6439}
  \end{pmatrix}$, \\

  characteristic polynomial=
  $101633401431659687926336+
        3020312682800x+x^2$.
%%%%%%%%%%%%%%%%%%%%%%%%%%%%%%%%%%%%%%%%%

\end{document}